\documentclass[reqno,11pt]{amsart}
\usepackage{graphicx, amsmath, amssymb, amscd, mathtools, hyperref, amsthm, euscript, amsfonts,bm,color}
\hypersetup{hypertexnames=false}
\usepackage{graphicx}
 
\usepackage{mathtools}
\usepackage{amsmath}
\usepackage{amssymb}
 
\usepackage{amsthm}
\usepackage{graphicx}
\usepackage{mathrsfs}
\usepackage{xcolor}

\usepackage{tikz-cd}
\usepackage{tikz}
\usetikzlibrary{patterns}
\usepackage{hyperref}
\usepackage{cleveref}
 
\usepackage[utf8]{inputenc}
\usepackage[T1]{fontenc}
\allowdisplaybreaks
\usepackage{verbatim}
\usepackage{float}
 
\usepackage{enumitem}
\newlist{condenum}{enumerate}{1}
\setlist[condenum]{label=(R\arabic*), ref=(R\arabic*), leftmargin=*}

\DeclareMathAlphabet{\mathpzc}{OT1}{pzc}{m}{it}
 
\usepackage{thmtools}
\usepackage{thm-restate}
 
\usepackage{caption}
 
\DeclareMathAlphabet{\mathpzc}{OT1}{pzc}{m}{it}
 
\numberwithin{equation}{section}
 
\newtheorem{theorem}{Theorem}[section]
 
\newtheorem{corollary}[theorem]{Corollary}

\newtheorem{cor}[theorem]{Corollary}

\newtheorem{proposition}[theorem]{Proposition}

\newtheorem{lemma}[theorem]{Lemma}
 
\theoremstyle{definition}

\newtheorem{Def}[theorem]{Definition}

\newcommand{\be}{begin{equation}}

\newcommand{\z}{\mathbb{Z}}

\newcommand{\br}{\mathbb{R}}

\newcommand{{\grinv}}{{\Cal G}^{-r}}

\newcommand{\ba}{\backslash}

\newcommand{\Cal}{\mathcal}

\newcommand{\SL}{\operatorname{SL}}

\newcommand{\bp}{\begin{pmatrix}}
\newcommand{\ep}{\end{pmatrix}}
\renewcommand{\be}{\begin{equation}}
\newcommand{\ee}{\end{equation}}
\renewcommand{\epsilon}{\varepsilon}

\renewcommand{\bp}{{\rm bp}}

\newcommand{\PS}{\op{PS}}\newcommand{\Haar}{\op{Haar}}

\newcommand{\norm}[1]{\lVert #1 \rVert}

\newcommand{\op}{\operatorname}
\newcommand{\supp}{\operatorname{supp}}

\newcommand{\BR}{\operatorname{BR}}
\newcommand{\BMS}{\operatorname{BMS}}

\renewcommand{\setminus}{-}

\newcommand{\Z}{\z}

\newcommand{\ga}{\gamma}
\newcommand{\F}{\mathcal F}
\newcommand{\La}{\Lambda}
\renewcommand{\i}{\op{i}}

\def\scrB{{\mathcal B}}

\def\scrL{{\mathcal L}}

\def\i{\mathrm{i}}

\def\dim{\operatorname{dim}}

\def\SL{\operatorname{SL}}

\def\supp{\operatorname{supp}}

\newcommand{\fa}{\mathfrak{a}}

\newcommand{\Ga}{\Gamma}\newcommand{\bb}{\mathbb}
\newcommand{\cal}{\mathcal}\newcommand{\NN}{\mathbb N}

\newcommand{\inte}{\op{int}}

\def\NN{{\mathbb N}}

\def\scrB{{\mathcal B}}

\def\scrL{{\mathcal L}}

\def\i{\mathrm{i}}

\def\dim{\operatorname{dim}}

\def\SL{\operatorname{SL}}

\def\supp{\operatorname{supp}}

\def\supp{\operatorname{supp}}

\newcommand{\al}{\alpha}
\newcommand{\te}{\theta}

\newcommand{\err}{\mathsf{E}}

\newcommand{\gap}{\mathsf{m}_\theta}
 
\begin{document}

\title[Haar local mixing for transverse subgroups]{Haar local mixing for transverse subgroups}

\author{Sam Edwards}
\address{Department of Mathematics, Uppsala University, Box 524, 751 20 Uppsala, Sweden}
\email{samuel.edwards@math.uu.se}
 
\author{Minju Lee}
\address{Department of Mathematical Sciences, KAIST, 291 Daehak-ro, Yuseong-gu, Daejeon 34141, Republic of Korea}
\email{minju.lee@kaist.ac.kr}
 
\author{Hee Oh}
\address{Department of Mathematics, Yale University, PO Box 208283, New Haven, CT 06520-8283, USA}
\email{hee.oh@yale.edu}
\thanks{Oh is partially supported by NSF grant DMS-2450703.}

\begin{abstract} 
Transverse subgroups form a special class of discrete subgroups of connected semisimple real algebraic groups that includes all discrete subgroups of rank-one Lie groups, as well as Anosov and relatively Anosov subgroups.
For torsion-free Zariski dense $\theta$-transverse subgroups, we establish Haar local mixing for the partial Weyl chamber action under a Bowen-Margulis-Sullivan (BMS) local mixing hypothesis.
The proof is based on a transfer principle using the domains of discontinuity introduced by Delarue-Monclair-Sanders and extended to transverse subgroups by Guedes Bonthonneau-Lefeuvre-Weich.
For $\theta$-Anosov subgroups, we verify the BMS hypothesis by adapting existing arguments and hence obtain Haar local mixing.
As a consequence, we obtain asymptotics for the volume and distribution of long returns and recover the growth indicator from return volumes.
Under linear $\theta$-regularity, the transfer also holds with polynomial error estimates.
This condition is satisfied by subgroups of $\theta$-Anosov or relatively $\theta$-Anosov groups.
\end{abstract}

\maketitle
\tableofcontents

\section{Introduction}

Local mixing describes the asymptotic behavior of suitably normalized correlations in infinite volume. 
For a discrete subgroup, Bowen-Margulis-Sullivan (BMS) measures reflect the geometry of its limit sets, while Haar measure describes volume in the ambient space.
A transfer from BMS local mixing to Haar local mixing therefore relates BMS dynamics to the volume and distribution of long returns.

In \cite{ELO}, we established this transfer for the full Weyl chamber action of Zariski dense discrete subgroups, adapting Roblin's rank-one approach \cite{Ro}. 
Here we extend it to partial Weyl chamber actions of torsion-free Zariski dense transverse subgroups. 
The main geometric difficulty is that the subgroup need not act properly on the corresponding homogeneous space. 
We work on the domains of discontinuity constructed by Delarue-Monclair-Sanders \cite{DMS1,DMS2} for Anosov subgroups and extended to transverse subgroups by Guedes Bonthonneau-Lefeuvre-Weich \cite[Section 7]{BLW}.
Their quotients provide smooth manifolds on which both measures and the partial Weyl chamber action are defined.

\subsection*{The partial Weyl chamber action}

To state our results, let $G$ be a connected semisimple real algebraic group. Fix a Cartan decomposition $G=K\exp(\mathfrak a^+)K$, let $\Pi$ be the corresponding set of simple roots, and fix a nonempty subset $\theta\subset\Pi$.
 
Throughout the paper, we assume that
$$\text{$\Gamma<G$ is torsion-free, Zariski dense and $\theta$-transverse. }$$
The definition is recalled in Section \ref{ps}.
In rank one, every discrete subgroup is $\Pi$-transverse. Every subgroup of a $\theta$-Anosov or relatively $\theta$-Anosov group is also $\theta$-transverse. 
The subset $\theta$ specifies the flow directions and the flag varieties used to describe the action.

Write
$$
 \mathfrak a_\theta=\bigcap_{\alpha\in\Pi\setminus\theta}\ker\alpha,\qquad \mathfrak a_\theta^+=\mathfrak a_\theta\cap\mathfrak a^+,\qquad A_\theta=\exp(\mathfrak a_\theta),
$$
and set $A_\theta^+=\exp(\mathfrak a_\theta^+)$.
Let $P_\theta=M_\theta A_\theta N_\theta^-$ be the corresponding standard parabolic subgroup with its Langlands decomposition, so that $Z_G(A_\theta)=M_\theta A_\theta$.
Let $\i:\mathfrak a\to\mathfrak a$ denote the opposition involution, and let $\Lambda_\theta$ and $\Lambda_{\i\theta}$ be the limit sets of $\Gamma$ in the flag varieties $\mathcal F_\theta=G/P_\theta$ and $\mathcal F_{\i\theta}=G/P_{\i\theta}$, respectively.
Let $\mathcal F_\theta^{(2)}$ denote the open $G$-orbit in $\mathcal F_\theta\times\mathcal F_{\i\theta}$, and set
$$
 \Lambda_\theta^{(2)}:=(\Lambda_\theta\times\Lambda_{\i\theta})\cap\mathcal F_\theta^{(2)}.
$$

The elements of $\mathcal F_\theta^{(2)}$ are pairs of flags in general position; in rank one, they are pairs of distinct points at infinity.
The Hopf parametrization identifies $G/M_\theta$ with $\mathcal F_\theta^{(2)}\times\mathfrak a_\theta$.
In these coordinates, right translation by $\exp v\in A_\theta$ adds $v$ to the $\mathfrak a_\theta$ coordinate.

When $\theta=\Pi$, the group $M_\theta$ is compact, so $\Gamma$ acts properly discontinuously on $G/M_\theta$ and the right $A_\theta$-action descends to the smooth quotient $\Gamma\backslash G/M_\theta$. For general $\theta$, however, $M_\theta$ can be noncompact, and $\Gamma$ need not act properly on $G/M_\theta$. Kim-Oh-Wang \cite[Theorem 9.1]{KOW} proved that $\Gamma$ acts properly discontinuously on $\widetilde\Omega=\Lambda_\theta^{(2)}\times\mathfrak a_\theta$. The quotient $\Omega=\Gamma\backslash\widetilde\Omega$ carries the Bowen-Margulis-Sullivan measures constructed from Patterson-Sullivan measures on the limit sets.

The space $\Omega$ supports the BMS dynamics. 
To compare these dynamics with ambient volume, we need a smooth neighborhood on which the same action is defined.
By \cite{DMS1,DMS2,BLW}, we may choose a $\Gamma$-invariant and $A_\theta$-invariant open neighborhood $\widetilde{\mathcal M}\subset G/M_\theta$ of $\widetilde\Omega$
on which $\Gamma$ acts properly discontinuously.
Since $\Gamma$ is torsion-free, this action is free, and
$$
 \mathcal M=\Gamma\backslash\widetilde{\mathcal M}
$$
is a smooth manifold containing $\Omega$. 
It carries the BMS measures, supported on $\Omega$, and Haar measure $dx$, induced by the $G$-invariant measure on $G/M_\theta$. 
The induced right $A_\theta$-action on $\mathcal M$ is called the \emph{partial Weyl chamber action}.
It preserves $\Omega$, and its restriction to $\Omega$ is the action underlying BMS local mixing.
We recall the required properties of this domain in Theorem \ref{theorem.dd}.

\subsection*{From BMS to Haar local mixing}

The BMS measure used in the mixing hypothesis is adapted to the direction of the flow. 
Let $\psi_\Gamma^\theta$ be the growth indicator of $\Gamma$, which records exponential growth in different Cartan directions.
Fix $u\in\mathfrak a_\theta^+\smallsetminus\{0\}$ and a $(\Gamma,\theta)$-proper linear form $\psi\in\mathfrak a_\theta^*$ tangent to $\psi_\Gamma^\theta$ at $u$ (see Section \ref{ps}).
We use this tangent form as the parameter for the Patterson-Sullivan measures defining our BMS measure.

Let $\nu_\psi$ and $\nu_{\psi\circ\i}$ be $(\Gamma,\psi)$- and $(\Gamma,\psi\circ\i)$-Patterson-Sullivan measures on $\Lambda_\theta$ and $\Lambda_{\i\theta}$, respectively. 
Denote the associated BMS, Burger-Roblin (BR) and dual Burger-Roblin measures on $\mathcal M$ by
$$
 m_\psi^{\BMS}=m^{\BMS}_{\nu_\psi,\nu_{\psi\circ\i}},\qquad m_\psi^{\BR}=m^{\BR}_{\nu_{\psi\circ\i}},\qquad m_\psi^{\BR_*}=m^{\BR_*}_{\nu_\psi}.
$$
Their definitions are given in Section \ref{sec.Norb}. If $\dim A_\theta>1$, all three measures have infinite total mass.
For $v\in\mathfrak a_\theta$, write $a_v:=\exp v$. 
We allow perturbations $\sqrt t\,v$, with $v\in\ker\psi$, of the translation vector $tu$.

Let $\Psi:(0,\infty)\to(0,\infty)$ be a normalizing function and $J:\ker\psi\to(0,\infty)$ a bounded Borel function. We say that $m_\psi^{\BMS}$ satisfies local mixing with the pair $(\Psi,J)$ if, for every $v\in\ker\psi$ and $f_1,f_2\in C_c(\mathcal M)$,
$$
 \lim_{t\to\infty} \Psi (t)\int_{\mathcal M}f_1(xa_{tu+\sqrt t\,v})f_2(x)\,dm_\psi^{\BMS}(x)=J(v)m_\psi^{\BMS}(f_1)m_\psi^{\BMS}(f_2).
$$

The function $J$ describes the dependence of the limit on $v$.
For fixed $v$, the normalized vectors $t^{-1}(tu+\sqrt t\,v)$ converge to $u$.

The case $v=0$ gives local mixing along the directional flow $x\mapsto xa_{tu}$. Since $J(0)>0$, this implies that the flow is topologically mixing on $\operatorname{supp}m_\psi^{\BMS}$, and hence topologically transitive.
It follows from \cite[Lemma 3.2]{ELO2} that $u\in\operatorname{int}\mathcal L_\theta$, where $\scrL_{\theta}$ is the limit cone of $\Gamma$, and the interior
is taken in $\mathfrak a_\theta$.

Passing to Haar measure introduces an exponential normalization.
Define
$$
 2\rho_\theta(H)
 :=\operatorname{tr}\bigl(\operatorname{ad}(H)|_{\mathfrak n_\theta^-}\bigr),
$$
where $\mathfrak n_\theta^-$ is the Lie algebra of $N_\theta^-$.
This linear form gives the logarithmic volume expansion of conjugation by $\exp H$ on $N_\theta^-$.
The difference $2\rho_\theta-\psi$ gives the additional exponential factor in the transfer theorem.

\begin{theorem}[Transfer from BMS to Haar local mixing]\label{thm.haar}
Suppose that $m_\psi^{\BMS}$ satisfies local mixing with the pair $(\Psi,J)$.
Then, for every  $v\in\ker\psi$ and $f_1,f_2\in C_c(\mathcal M)$,
$$\lim_{t\to\infty}\Psi(t)e^{(2\rho_\theta-\psi)({tu+\sqrt t\,v})}\int_{\mathcal M}f_1(xa_{tu+\sqrt t\,v})f_2(x)\,dx=J(v)m_\psi^{\BR}(f_1)m_\psi^{\BR_*}(f_2).$$
\end{theorem}
The same functions $\Psi$ and $J$ appear in the BMS and Haar local mixing limits.
The two limiting factors are now the BR and dual BR measures: $m_\psi^{\BR}$ describes the final position, and $m_\psi^{\BR_*}$ the initial position, in the Haar correlation.

\subsection*{Long returns and the growth indicator}

Haar local mixing gives a geometric description of long returns.
Write $\operatorname{vol}(E)=\int_E dx$, and let $B(p,r)\subset\mathcal M$ be the ball of radius $r$ centered at $p$ for any fixed smooth Riemannian metric on $\mathcal M$.
The set $B(p,r)\cap B(q,s)a_{-tu}$ consists of points in $B(p,r)$ whose images at time $t$ lie in $B(q,s)$.
The following corollary determines its volume and recovers the directional growth indicator from its exponential rate.

\begin{corollary}\label{cor.vol}
 Suppose that $m_\psi^{\BMS}$ satisfies local mixing with the pair $(\Psi,J)$.
The following hold for any $p,q\in\operatorname{supp}m_\psi^{\BMS}$ and sufficiently small $r,s>0$.
 \begin{enumerate}
     \item 
 As $t\to\infty$,
\begin{multline}\label{eq.rv}
 \operatorname{vol}\bigl(B(p,r)\cap B(q,s)a_{-tu}\bigr)\sim  \\ {J(0)m_\psi^{\BR_*}(B(p,r))m_\psi^{\BR}(B(q,s))}{\Psi(t)}^{-1}e^{t(\psi_\Gamma^\theta(u)-2\rho_\theta(u))}.
\end{multline}
 
\item We have 
\begin{equation}\label{eq.growth-volume} \psi_\Gamma^\theta(u)=2\rho_\theta(u)+\lim_{t\to\infty}\frac{1}{t}\log\operatorname{vol}\bigl(B(p,r)\cap B(q,s)a_{-tu}\bigr). \end{equation}
 \end{enumerate}
 
\end{corollary}

Tangency gives $\psi(u)=\psi_\Gamma^\theta(u)$.
Moreover, local mixing implies $\log\Psi(t)=o(t)$ (Lemma \ref{lem.ratio}), so the normalization does not contribute to the exponential rate in \eqref{eq.growth-volume}.

\subsection*{Polynomial error estimates}

The transfer theorem requires no linear regularity assumption.
For polynomial error estimates, we need a quantitative form of the
contraction used in its proof. Let $\mu$ and $\mu_\theta$ denote the
full and partial Cartan projections (see Section \ref{ps}).
We say that $\Gamma$ is \emph{linearly $\theta$-regular} if there
exists $C>1$ such that, for every $\gamma\in\Gamma$,

\begin{equation}\label{eq.rc}
 \min_{\alpha\in\theta}\alpha(\mu(\gamma))\geq C^{-1}\|\mu_\theta(\gamma)\|-C.
\end{equation}

Linear $\theta$-regularity holds if $\Gamma$ is contained in a $\theta$-Anosov or relatively $\theta$-Anosov subgroup; see Section \ref{ps}.

\begin{theorem}[Polynomial error transfer] \label{thm.eh}
Suppose that $\Gamma$ is linearly $\theta$-regular and that $m_\psi^{\BMS}$ satisfies local mixing with the pair $(\Psi,J)$ and with error term
$$
 O_{\mathcal K,v}\!\left(t^{-d}\|f_1\|_{C^k(\mathcal M)}\|f_2\|_{C^k(\mathcal M)}\right)
$$
for some $d>0$ and integer $k\geq1$, every compact $\mathcal K\subset\mathcal M$, and every $v\in\ker\psi$.
Then there exists $d'>0$ such that, for every such $\mathcal K$ and $v$, the error in \Cref{thm.haar} is
$$
 O_{\mathcal K,v}\!\left(t^{-d'}\|f_1\|_{C^{k}(\mathcal M)}\|f_2\|_{C^{k}(\mathcal M)}\right).
$$
Both estimates hold as $t\to\infty$, uniformly for $f_1,f_2\in C_c^\infty(\mathcal K)$.
\end{theorem}

\subsection*{Applications}

The transfer theorem applies whenever the corresponding BMS
local mixing is known. We recall several such cases.

For geometrically finite groups in rank one, the BMS measure is finite (\cite{Sul}, \cite{CI}).
Babillot \cite{Bab} proved mixing for finite BMS measures, so the local mixing hypothesis holds with $\Psi(t)=1$.
For infinite-index co-abelian subgroups of convex cocompact groups, the BMS measure is infinite, and Oh-Pan \cite{OP} proved local mixing with $\Psi(t)=t^{d/2}$, where $d$ is the rank of the abelian quotient. 
In higher rank, Thirion \cite{ThiMix} proved local mixing for Ping-Pong groups along the maximal growth direction.
Sambarino \cite{Samb1} established BMS local mixing for Zariski dense hyperconvex representations of fundamental groups of closed negatively curved manifolds. His later work \cite{Samb4} gives directional local mixing for general $\theta$-Anosov groups.
Chow-Sarkar \cite{CS} proved the BMS local mixing hypothesis for Borel Anosov groups.
In Appendix \ref{sec.app}, we show that their proof also applies to general $\theta$. More precisely, if $\Gamma$ is $\theta$-Anosov and $u\in\operatorname{int}\mathcal L_\theta$, then the BMS local mixing hypothesis holds with
$$
\Psi(t)=t^{(\dim A_\theta-1)/2},\qquad J(v)=\kappa_u e^{-Q_u(v)},
$$
where $\kappa_u>0$ and $Q_u$ is a positive definite quadratic form on $\ker\psi$.
 
Combining this with Theorem \ref{thm.haar} gives the following, which applies in particular to the images of Zariski
dense Hitchin representations of closed surface groups for every nonempty $\theta\subset\Pi$ \cite{Lab}.
\begin{theorem}[Haar local mixing for $\theta$-Anosov subgroups]
Let $\Gamma<G$ be a torsion-free Zariski dense $\theta$-Anosov subgroup and $u\in \inte\mathcal L_\theta$.
There exist $\kappa_u>0$ and a positive definite quadratic form $Q_u$ on $\ker\psi$ such that, for every $v\in\ker\psi$ and $f_1,f_2\in C_c(\mathcal M)$,
\begin{multline*}
 \lim_{t\to\infty}t^{(\dim A_\theta-1)/2}e^{(2\rho_\theta-\psi)(tu+\sqrt t\,v)}\int_{\mathcal M}f_1(xa_{tu+\sqrt t\,v})f_2(x)\,dx \\ 
 =\kappa_u e^{-Q_u(v)}m_\psi^{\BR}(f_1)m_\psi^{\BR_*}(f_2).
\end{multline*}
\end{theorem}

If $\Gamma$ is $\theta$-Anosov and $u\in\operatorname{int}\mathcal L_\theta$, then, for $p,q,r,s$ as in Corollary \ref{cor.vol},

\begin{equation}\label{avol}
\operatorname{vol}\bigl(B(p,r)\cap B(q,s)a_{-tu}\bigr)\sim C\,t^{-(\dim A_\theta-1)/2} e^{-t(2\rho_\theta(u)-\psi_\Gamma^\theta(u))}
\end{equation}
as $t\to\infty$, for some $C=C(u,p,q,r,s)>0$.

For Anosov subgroups, the BR measures also have a spectral interpretation.
Guedes Bonthonneau-Lefeuvre-Weich \cite[Section 4.3]{BLW} realize them as leading resonant and coresonant states, represented by distributional differential forms, whose pairing gives the BMS measure.
This relates the two BR factors in our Haar local mixing formula to the spectral theory of the partial Weyl chamber action.

For Zariski dense Anosov subgroups satisfying a Diophantine condition, Guedes Bonthonneau-Lefeuvre-Weich \cite[Theorem 1.8]{BLW} obtained a full asymptotic expansion for BMS correlations along a fixed direction. 
Kim-Oh-Pan \cite{KOP} established BMS local mixing with polynomial error estimates for several classes of relatively Anosov groups, including images of cusped positive representations.
Theorem \ref{thm.eh} therefore gives corresponding polynomial error estimates for Haar local mixing.

\subsection*{Proof strategy and organization}

The proof proceeds through two horospherical equidistribution statements, following \cite[Sections 5-6]{ELO} (see also \cite{OhShah}).
BMS local mixing first gives equidistribution of translates of Patterson-Sullivan measures on horospherical leaves. 
Comparing these with translates of Lebesgue measures gives equidistribution with respect to the BR measure.
Integration then gives Haar local mixing.

The transfer argument passes from integrals on horospherical leaves to integrals over small neighborhoods of those leaves. 
For fixed $v\in\ker\psi$ and $f\in C_c(\mathcal M)$, this requires controlling
\begin{equation}\label{con}
 \bigl|f(\Gamma ph a_{tu+\sqrt t\,v}M_\theta)-f(\Gamma p a_{tu+\sqrt t\,v}M_\theta)\bigr|,
\end{equation}
where $p$ ranges over a fixed compact subset of $G$ projecting into $\widetilde{\mathcal M}$, and $h$ lies in a fixed small compact neighborhood of the identity in $N_\theta^-$.
The main additional difficulty is the noncompactness of $M_\theta$ when $\theta\ne\Pi$. At a return to a fixed compact subset of $\mathcal M$, we can write
$$
 \gamma p a_{tu+\sqrt t\,v}=qm,\qquad \gamma\in\Gamma,\quad m\in M_\theta,
$$
with $q$ in a fixed compact subset of $G$. Then
$$
 \gamma ph a_{tu+\sqrt t\,v}M_\theta=q\bigl(m a_{tu+\sqrt t\,v}^{-1}h a_{tu+\sqrt t\,v}m^{-1}\bigr)M_\theta.
$$
When $\theta=\Pi$, the factor $m$ belongs to a compact group. 
For general $\theta$, it need not remain bounded, so contraction by $a_{tu+\sqrt t\,v}^{-1}$ alone does not suffice.
 
Lemma \ref{lem.return} shows that $m a_{tu+\sqrt t\,v}^{-1}h a_{tu+\sqrt t\,v}m^{-1}$ tends to the identity uniformly over all such returns and $h$. 
By continuity of $f$ and compactness of the set of representatives $q$, the difference in \eqref{con} therefore tends to zero. Under linear $\theta$-regularity, the conjugates approach the identity exponentially in $t$.
Together with the smoothing estimates, this proves the polynomial error estimate in Theorem \ref{thm.eh}.
 
\medskip 
Section \ref{ps} recalls transverse subgroups and PS measures. In Section \ref{sec.meas} we define $\widetilde{\mathcal M}$ and prove Lemma \ref{lem.return}.
 Section \ref{sec.Norb} defines the generalized BMS measures and gives the product formulas for these measures. Section \ref{sec.lcd} contains the smoothing estimates. Sections \ref{sec.set} and \ref{s2} prove the two horospherical equidistribution statements and the Haar transfer theorem.

\subsection*{AI statement}
We used AI tools to assist with proofreading, reference searches, and  checking proofs. 
\section{Transverse subgroups and Patterson-Sullivan measures}\label{ps}
Let $G$ be a connected semisimple real algebraic group, and let $\Gamma<G$ be a Zariski dense discrete subgroup.
Fix a Cartan decomposition $\mathfrak g=\mathfrak k\oplus\mathfrak p$ of the Lie algebra of $G$ \cite[Chapter VI (6.23)]{Knapp}, and let $K$ be the maximal compact subgroup with Lie algebra $\mathfrak k$. Choose a maximal abelian subspace $\mathfrak a$ of $\mathfrak p$.
Fix a left $G$-invariant and right $K$-invariant Riemannian metric on $G$, and denote its distance by $d$. The inner product at the identity restricts to a Weyl-group invariant inner product on $\mathfrak a$, denoted by $\langle\cdot,\cdot\rangle$, with corresponding norm $\|\cdot\|$.
 
Let $A:=\exp \mathfrak a$. Choosing a closed positive Weyl chamber $\fa^+$ of $\fa$, let $A^+=\exp\mathfrak a^+$. 
Denote by $\Pi$ the set of simple roots corresponding to $\fa^+$. Fix a nonempty subset $\theta\subset \Pi$, and let
\begin{equation*}
 \fa_\theta:=\bigcap _{\al\in \Pi\setminus\theta}\ker \al,\;\; \;\; \fa_\theta^+:=\fa_\theta\cap \fa^+,\;\; \;\; A_\theta:=\exp\fa_\theta,\;\; \;\; A_\theta^+=\exp \mathfrak
a_\theta^+.
\end{equation*}
We write $\op{int}A_\te^+$ and $\op{int}\mathfrak a_\te^+$ for the interiors relative to $A_\te$ and $\mathfrak a_\te$, respectively.
For any $a\in\operatorname{int}A_\theta^+$, define
$$
 N_\theta^\pm=\{g\in G:a^{-n}ga^n\to e\text{ as }n\to\mp\infty\}.
$$
These groups are independent of the choice of $a$. We write $N_\theta:=N_\theta^-$. Denote their Lie algebras by $\mathfrak n_\theta^\pm$.
Let $L_\theta=Z_G(A_\theta)$, and let $M_\theta$ be the subgroup generated by $K\cap L_\theta$ and the derived subgroup $[L_\theta, L_\theta]$. Then $L_\theta=A_\theta
M_\theta$. 
We set
$$ P_\theta:=P_\te^-=L_\theta N_\theta^-,\quad 
 P_\theta^+:=L_\theta N_\theta^+, \quad \F_\theta:=G/P_\theta.
$$
These are opposite parabolic subgroups with common Levi subgroup $L_\theta$, and $P_\theta=M_\theta A_\theta N_\theta^-$ is its Langlands decomposition \cite[Chapter VII. $\mathsection$7]{Knapp}.
Since $G=KP_\theta$ and $K\cap P_\theta=K\cap M_\theta$, the map $k(K\cap M_\theta)\mapsto kP_\theta$ induces a diffeomorphism $K/(K\cap M_\theta)\simeq\mathcal F_\theta$ \cite[Proposition 7.83]{Knapp}.
 
Let $w_0\in N_K(A)$ represent the longest Weyl element. The opposition involution is $\i=-\operatorname{Ad}(w_0)|_{\fa}$; it preserves $\fa^+$ and induces a permutation of $\Pi$. 
We write $\i\theta$ for the image of $\theta$.  
A pair $(\xi, \eta)\in \F_\theta \times \F_{\i\theta}$ is said to be in general position if
it lies in the unique open $G$-orbit.
Let $p_\theta:\mathfrak a\to\mathfrak a_\theta$ be the orthogonal projection with respect to $\langle\cdot,\cdot\rangle$.
 
\begin{Def}
The \textit{partial Cartan projection} $\mu_\te:G\to \fa_\te^+$ is defined as
$$
 \mu_\te(g)=p_\te (\mu(g))
$$
for each $g\in G$, where $\mu(g)\in \mathfrak a^+$ is the unique element such that $g\in
K\exp(\mu(g))K$ \cite[Theorem 7.39]{Knapp}.     
\end{Def}
The {\it $\mathfrak a$-valued Busemann map} $\beta: \cal F_\Pi \times G \times G \to\mathfrak a $ is defined as
$$
 \beta_\xi ( g, h):=\sigma (g^{-1}, \xi)-\sigma(h^{-1}, \xi)
$$
for $\xi\in \cal F_\Pi$ and $g, h\in G$;
here $\sigma(g^{-1},\xi)\in \fa$ is the unique element such that $g^{-1}k \in K \allowbreak\exp
(\sigma(g^{-1}, \xi)) N_{\Pi}$ if $\xi=kP_{\Pi}$ for some (and hence all) $k\in K$. 
\begin{Def}
For $(\xi,g,h)\in \cal F_\theta\times
G\times G$, we define
$$
 \beta_{\xi}^\theta (g, h): = p_\theta ( \beta_{\xi_0} (g, h))
$$
for any $\xi_0\in \F_\Pi$ projecting to $\xi$. This is independent of the choice of $\xi_0$
\cite[Lemma 6.1]{Quint2}.    
\end{Def}
We also write $\sigma_\theta(g, \xi)=p_\theta (\sigma(g, \xi_0))$. 
Since both the natural projection $\mathcal F_\Pi\to\mathcal F_\theta$ and $\sigma$ are smooth \cite[Theorem 6.46]{Knapp}, the Busemann maps are smooth.
It is straightforward to check that the following identities hold \cite[Lemma 6.2]{Quint2}:
\begin{equation}\label{eq.basic0}
\beta_\xi^\te(g,h)+\beta_\xi^\te(h,q)=\beta_\xi^\te(g,q),\quad \beta_{g\xi}^\te(gh,gq)=\beta_\xi^\te(h,q).
\end{equation}
 
\begin{Def}
The limit set $\Lambda_\theta=\Lambda_\theta(\Gamma)$ consists of all $\xi\in\mathcal F_\theta$ for which there is a sequence $\gamma_n\to\infty$ in $\Gamma$, with Cartan decompositions $\gamma_n=k_n\exp(\mu(\gamma_n))\ell_n$, $k_n,\ell_n\in K$, such that $\min_{\alpha\in\theta}\alpha(\mu(\gamma_n))\to+\infty$ and $k_nP_\theta\to\xi$ in $\mathcal F_\theta$.
\end{Def}
This definition agrees with that of Benoist \cite{Ben}.

\begin{Def} \label{def.conf}\rm
Given $\psi\in \mathfrak a_\theta^*$ and a subgroup $\Delta <G$, a Borel probability measure $\nu_{\psi}$ on $\F_\theta$ is called a $(\Delta, \psi)$-conformal measure if, for every $\ga\in \Delta$ and $\nu_\psi$ almost every $\xi\in \F_\theta$,
\be\label{gc0}
 \frac{d\ga_* \nu_\psi}{d\nu_\psi}(\xi) =e^{\psi (\beta_\xi^\theta (e, \gamma))},
\ee
where $\ga_* \nu_\psi(Q)=\nu_{\psi}(\gamma^{-1} Q)$ for any Borel subset $Q\subset \F_\theta$. When $\Delta=\Gamma$ and it is supported on $\La_\theta$, we call it a $(\Ga, \psi)$-PS measure.
\end{Def}
Let $\rho_\te$ denote the half sum of all positive roots, counted with multiplicity, with respect to $\mathfrak a_{\te}^+$; i.e., for all $H\in \fa_\theta$,
$$
2\rho_\theta(H)=\op{tr}(\op{ad}H|_{\mathfrak n_\theta^-}).
$$
By \cite[Lemma 6.3]{Quint2}, the $K$-invariant probability measure $m_\te$ on $\F_\theta$ is a $(G, 2\rho_\te)$-conformal measure.

\begin{Def}\label{def.lc}
The limit cone $\cal L_\te=\cal L_\te(\Ga)$ is defined by
$$
 \cal L_\te := \left\{ \lim_{n\to\infty}t_n\mu_\te(\gamma_n): \gamma_n\in\Ga,\ t_n>0,\ t_n\to0 \right\}.
$$    
\end{Def}
The cone $\cal L_\te$ has nonempty interior in $\fa_\te$ \cite[Theorem 1.2]{Ben}.
 
\begin{Def}[Transverse group]\rm
A discrete subgroup $\Ga<G$ is called 
\begin{enumerate}
    \item \textit{$\te$-regular} if $\min_{\alpha\in\te}\alpha(\mu(\gamma))\to+\infty$  as $\gamma\to\infty$ in $\Ga$, and 
    \item \textit{$\te$-antipodal} if any two distinct points of $\La_{\te\cup\i\te}$ are in general position.
\end{enumerate}
We say $\Ga$ is \emph{$\te$-transverse} if it is both $\te$-regular and $\te$-antipodal.
\end{Def}
In particular, if $\Ga$ is $\te$-transverse, then $\mu_\te|_\Ga$ is proper.
 
\begin{Def}[Growth indicator] \rm
Let $\Ga<G$ be a $\te$-transverse subgroup. 
For $u\in\fa_\te\setminus\{0\}$,
\begin{equation*}
 \psi_\Ga^\te(u) := \|u\|
 \inf_{\substack{\cal C\subset\fa_\te\text{ open cone},\; u\in\cal C}}\tau_{\cal C},
\end{equation*}
where $\tau_{\cal C}$ is the abscissa of convergence of the series $ \sum_{\substack{\gamma\in\Ga,
\mu_\te(\gamma)\in\cal C}} e^{-s\|\mu_\te(\gamma)\|}$.  
\end{Def}
 
We extend $\psi_\Ga^\te$ to all of $\mathfrak a_\te$ by defining $\psi_\Ga^\te(0):=0$. We have
$
\psi_\Ga^\te\le 2\rho_\te
$ by \cite[Theorem 1.4(1)]{KOW}, applied to $m_\te$.
 
\begin{theorem}[{\cite[Theorem 1.1.1]{Quint1}, \cite[Theorem 3.3]{KOW}}]
Let $\Ga<G$ be a Zariski dense $\te$-transverse subgroup. The function $\psi_\Ga^\te:\fa_\te\to\bb
R\cup\{-\infty\}$ is homogeneous of degree one, upper semicontinuous, and concave. Moreover,
$$
 \cal L_\te =\{u\in\fa_\te:\psi_\Ga^\te(u)\geq0\} =\{u\in\fa_\te:\psi_\Ga^\te(u)>-\infty\},
$$
and $\psi_\Ga^\te$ is positive on $\op{int}\cal L_\te$.
\end{theorem}

A linear form $\psi\in\fa_\te^*$ is called
\begin{enumerate}
    \item \textit{$(\Ga,\te)$-proper} if  
$\psi(\mu_\te(\gamma))\to+\infty$ as $\gamma\to \infty$ in $\Gamma$.
    \item \textit{a tangent form} if $\psi\ge
\psi_\Ga^\theta$ and $\psi(v)=\psi_\Ga^\theta(v)$ for some nonzero $v\in \cal L_\theta$; in this case we say that $\psi$ is tangent to $\psi_\Ga^\theta$ at $v$.
\end{enumerate}

If a linear form $\psi$ is positive on $\cal L_\theta-\{0\}$, then it is $(\Ga, \theta)$-proper \cite[Lemma 4.3(1)]{KOW}.

\begin{theorem}[{\cite[Theorem 8.4]{Quint2}, \cite[Proposition 5.12]{KOW}}]\label{thm.PS}
Suppose that $\Ga$ is $\te$-transverse.
If $\psi$ is tangent to $\psi_\Ga^\te$ at some $v\in\op{int}\mathfrak a_\te^+$, then there exists a $(\Ga,\psi)$-PS measure.
\end{theorem}

\subsection*{Anosov and relatively Anosov subgroups}
Anosov representations were introduced by Labourie \cite{Lab} for surface groups and extended to word hyperbolic groups by Guichard-Wienhard \cite{GW}.
Anosov and relatively Anosov subgroups generalize, respectively, convex cocompact and geometrically finite subgroups of rank-one Lie groups.
Every Zariski dense subgroup of a $\theta$-Anosov or relatively $\theta$-Anosov subgroup is $\theta$-transverse.
 
A $\theta$-transverse group $\Gamma$ is called
\textit{relatively $\theta$-Anosov} if every
$\xi\in\Lambda_\theta$ satisfies one of the following:
\begin{enumerate}
\item $\xi$ is conical: there are distinct $a,b\in\Lambda_\theta$ and a sequence $\gamma_n\to\infty$ in $\Gamma$ such that $\gamma_n\xi\to a$ and $\gamma_n\zeta\to b$ uniformly for $\zeta$ in compact subsets of $\Lambda_\theta\setminus\{\xi\}$.
\item $\xi$ is bounded parabolic: its stabilizer is infinite and acts properly on $\Lambda_\theta\setminus\{\xi\}$ with compact quotient. We also require this stabilizer
to be finitely generated.
\end{enumerate}
See \cite[Definition~7.7, Theorem~7.10]{KL}.
The \textit{$\theta$-Anosov} case is the special case in which every point of $\Lambda_\theta$ satisfies (1) \cite[Section~3.3, Theorem~1.1]{KLP}.
 
In addition, every Zariski dense subgroup $\Gamma<G$ contained in a $\theta$-Anosov or relatively $\theta$-Anosov subgroup is linearly $\theta$-regular, in the sense of \eqref{eq.rc}.
Indeed, if $\Ga$ is itself $\te$-Anosov, then this follows from uniform regularity \cite[Theorem 1.1]{KLP}.
If $\Ga$ is relatively $\te$-Anosov, then it follows from the singular value estimates and the reduction in \cite[Theorem 1.7, Propositions 5.1, 13.4]{ZZ1}.
Clearly \eqref{eq.rc} then holds for every subgroup.
 
The assumption \eqref{eq.rc} is needed only for the polynomial error estimates.
 
\section{The domain of discontinuity and compact returns}\label{sec.meas}
In the rest of the paper, let $\Ga<G$ be a Zariski dense torsion-free $\theta$-transverse subgroup. Set $o:=M_\theta\in G/M_\theta$.
For $gM_\theta \in G/M_\theta$, define 
\begin{equation}\label{eq.hopf0}
 g^+:=gP_\te\in\F_\te, \qquad g^-:=gw_0P_{\i\te}\in\F_{\i\te}.
\end{equation}
Set $\F_\te^{(2)} :=G\cdot(e^+,e^-) =\{(g^+,g^-):g\in G\},$ which is the unique open $G$-orbit in $\F_\te \times \F_{\i\te}$ \cite[Section 3.3]{Ben}.
 
\begin{Def}\label{def.hopf}
\rm
The diffeomorphism $G/M_\te \to\F_\te^{(2)}\times\fa_\te$ given by
$$
 gM_\te \mapsto(g^+,g^-,\beta_{g^+}^{\te}(e,g))
$$
is called the \textit{Hopf parametrization}.
\end{Def}
For $g\in G$, $a\in A$, $m\in M_\theta$, $n\in N_\theta^+$ and $h\in N_\theta^-$,  \eqref{eq.hopf0} and \eqref{eq.basic0} give
\begin{equation}\label{eq.basic1}
 \begin{aligned}
 (gnam)^-&=g^-,& \beta_{g^-}^{\i\theta}(e,gnam) &=\beta_{g^-}^{\i\theta}(e,g)-\i p_\theta(\log a),\\
 (gham)^+&=g^+,& \beta_{g^+}^{\theta}(e,gham) &=\beta_{g^+}^{\theta}(e,g)+p_\theta(\log a).
 \end{aligned}
\end{equation}
In particular, for every $g\in G$ and $m\in M_\theta$, we have $(gm)^\pm=g^\pm$ and
$$
\beta_{g^+}^\theta(e,gm)=\beta_{g^+}^\theta(e,g),\qquad \beta_{g^-}^{\i\theta}(e,gm)=\beta_{g^-}^{\i\theta}(e,g).
$$
Thus both $\beta_{g^+}^\theta(e,g)$ and $\beta_{g^-}^{\i\theta}(e,g)$ depend only on the coset $gM_\theta$.
 
In the Hopf parametrization, the left $G$-action and the right $A_\te$-action are given, for $g\in G$ and $v\in\mathfrak a_\te$, by
\begin{equation}\label{eq.hopfa}
 g(\xi,\eta,b) =(g \xi,g \eta,b+\beta_{g \xi}^{\te}(e,g)),\qquad (\xi,\eta,b)\exp v=(\xi,\eta,b+v).
\end{equation}
 
\medskip
\noindent{\bf The domain of discontinuity.}
Set
$$
 \La_\te^{(2)} :=(\La_\te\times\La_{\i\te})\cap\F_\te^{(2)}, \qquad \widetilde\Omega:=\La_\te^{(2)}\times\fa_\te\subset G/M_\te.
$$
Since $\Ga$ is $\te$-transverse, the natural projections from $\La_{\te\cup\i\te}$ to $\La_\te$ and $\La_{\i\te}$ are $\Ga$-equivariant homeomorphisms \cite[Lemma 9.5]{KOW}. 
Write 
$$
\La_\te \to \La_{\i\te}\qquad \xi\mapsto \xi^*    
$$
for the resulting identification.
For
$\xi,\eta\in\La_\te$, the flags $\xi$ and $\eta^*$ are in general position if and only if
$\xi\ne\eta$.
 
Following the construction of Delarue-Monclair-Sanders \cite{DMS1,DMS2} and its extension
in \cite[Section 7]{BLW}, define the domain $\widetilde{\mathcal M}\subset G/M_\te$ by
\begin{equation}\label{eq.om}
 \widetilde{\mathcal M}:=\bigl\{gM_\te: \text{for every }\xi\in\La_\te, (g^+, \xi^*)\in\F_\te^{(2)} \ \text{or}\ (\xi,g^-)\in\F_\te^{(2)}\bigr\}.
\end{equation}
For $\ga\in\Ga$, set $\gap(\gamma):=\min_{\alpha\in\theta}\alpha(\mu(\gamma))$. 
Fix norms on $\mathfrak n_\theta^\pm$.
The following estimates also follow from \cite[Lemma 7.10]{BLW}.
We give a proof adapting the arguments of \cite[Proposition 9.7]{KOW} and \cite[Lemma 5.2]{SF}.
\begin{lemma}[Compact returns]\label{lem.return}
For every compact set $L_0\subset G$ with $L_0.o\subset\widetilde{\mathcal M}$, there exists $C=C(L_0)>0$ such that, for every $\gamma\in\Ga$, $v\in\fa_\theta$ and $m\in M_\theta$ satisfying $\gamma L_0\exp(v)\cap L_0m\neq\emptyset$,
\begin{equation}\label{eq.return1}
 \min \{ \|v-\mu_\theta(\gamma^{-1})\|, \|v+\mu_\theta(\gamma)\| \}\leq C.
\end{equation}
Moreover, if $v\in\fa_\theta^+$, then $\|v-\mu_\theta(\gamma^{-1})\|\leq C$, and
\begin{equation}\label{eq.holo}
 \max\left\{\left\|\operatorname{Ad}(m\exp(-v))\big|_{\mathfrak n_\theta^-}\right\|,\left\|\operatorname{Ad}(\exp(v)m^{-1})\big|_{\mathfrak n_\theta^+}\right\|\right\}\leq Ce^{-\gap(\gamma^{-1})}.
\end{equation}
\end{lemma}
\begin{proof}
Choose $p,q\in L_0$ such that $\gamma p\exp(v)=qm.$
Using the Cartan decomposition of $M_\theta$, write
$$
 m^{-1}=k_1\exp(b)k_2,\qquad k_1,k_2\in K\cap M_\theta,\quad b\in\ker p_\theta.
$$
Choose a Weyl representative $w$ from a fixed finite subset of $N_K(A)$ such that $v+b=\operatorname{Ad}(w)H$ with $H\in\fa^+$.
Then
$$
 \gamma^{-1}=pk_1w\exp(H)w^{-1}k_2q^{-1}
$$
and $\mu(\gamma^{-1})=H+O_{L_0}(1)$ by \cite[Lemma 4.6]{Ben}.
 
We claim that there is a finite set $F\subset\Ga$ such that $w\in M_\theta\cup M_\theta w_0$ whenever $\gamma\notin F$.
Suppose not.
Then there is a sequence such that $\ga_jp_j\exp(v_j)=q_jm_j$ with
$\gamma_j\to\infty$ and $w_j\notin M_\theta\cup M_\theta w_0$.
Passing to a subsequence, assume that $p_j\to p$, $q_j\to q$, $k_{1,j}\to k_1$, $k_{2,j}\to k_2$, and $w_j=w$.
 
Write
$$
\gamma_j^{-1}=k_j\exp(\mu(\gamma_j^{-1}))\ell_j,
\qquad k_j,\ell_j\in K.
$$
Since $H_j=\mu(\gamma_j^{-1})+O_{L_0}(1)$ and $\Gamma$ is $\theta\cup\i\theta$-regular, we have $\min_{\alpha\in\theta\cup\i\theta}\alpha(H_j)\to+\infty$.
As $p_jk_{1,j}w\to pk_1w$ and $w^{-1}k_{2,j}q_j^{-1}$ remains bounded, \cite[Lemma 2.5(1)-(2)]{SF}, applied to
$$
\gamma_j^{-1} =p_jk_{1,j}w\exp(H_j)w^{-1}k_{2,j}q_j^{-1},
$$
gives $k_jP_{\theta\cup\i\theta}\to pk_1wP_{\theta\cup\i\theta}$.
Hence $pk_1wP_{\theta\cup\i\theta} \in\Lambda_{\theta\cup\i\theta}$ by the definition of the limit set.
 
Since $pk_1.o=p.o\in\widetilde{\mathcal M}$, \eqref{eq.om} implies that either
$$
 (wP_\theta,w_0P_{\i\theta})\in\F_\theta^{(2)} \quad\text{or}\quad (P_\theta,wP_{\i\theta})\in\F_\theta^{(2)}.
$$
By \cite[Lemma 2.7]{SF}, applied also with $\theta$ replaced by $\i\theta$, these imply respectively $w\in M_\theta$ or $w\in M_\theta w_0$.
This is a contradiction and the claim follows.
 
For $\gamma\notin F$, projecting $v+b=\operatorname{Ad}(w)H$ onto $\fa_\theta$ yields
\begin{align}\label{eq.vph}
 v&=p_\theta(H)
   =\mu_\theta(\gamma^{-1})+O_{L_0}(1),
     \quad \text{if }w\in M_\theta,\notag\\
 v&=p_\theta(\operatorname{Ad}(w_0)H)
   =-\mu_\theta(\gamma)+O_{L_0}(1),
     \quad \text{if }w\in M_\theta w_0.    
\end{align}
This proves that \eqref{eq.return1} holds for all $\gamma\not\in F$.
In the second case when $w\in M_\te w_0$, $v\in-\fa_\theta^+$.
If also $v\in\fa_\theta^+$, then $v=0$, and  \eqref{eq.vph} bounds $\mu_\theta(\gamma)$.
Since $\mu_\theta|_\Ga$ is proper, there are only finitely many such $\gamma$.
Enlarging $F$, we may assume that $w\in M_\theta$ whenever $\gamma\not\in F$ and $v\in\fa_\theta^+$.
Therefore, whenever $\gamma\notin F$ and $v\in\fa_\theta^+$, $\|v-\mu_\theta(\gamma^{-1})\|\leq C$ by the first line of \eqref{eq.vph}.
 
We now prove the lemma for $\ga\not\in F$.
It remains to prove \eqref{eq.holo}.
Assume that $v\in\fa_\theta^+$.
Since $w\in M_\theta$, we have
$$
 \begin{aligned}
 m\exp(-v)
   &=k_2^{-1}w\exp(-H)w^{-1}k_1^{-1},\\
 \exp(v)m^{-1}
   &=k_1w\exp(H)w^{-1}k_2.
 \end{aligned}
$$
Note that the factors $k_1$, $k_2$, $w$ lie in a fixed compact subset of $M_\te$ and normalize $\mathfrak n_\theta^\pm$.
The root-space decomposition therefore gives
$$
 \begin{aligned}
 &\max\left\{
 \left\|\operatorname{Ad}(m\exp(-v))
       \big|_{\mathfrak n_\theta^-}\right\|,
 \left\|\operatorname{Ad}(\exp(v)m^{-1})
       \big|_{\mathfrak n_\theta^+}\right\|
 \right\}\\
 &\qquad\ll
 e^{-\min_{\alpha\in\theta}\alpha(H)}
 \ll_{L_0}e^{-\gap(\gamma^{-1})}.
 \end{aligned}
$$
Finally, suppose that $\gamma\in F$.
Since $\gamma^{-1}q=p\exp(v)m^{-1}$, \eqref{eq.basic1} implies that
$$
 v=\beta_{p^+}^\theta(e,\gamma^{-1}q)
   -\beta_{p^+}^\theta(e,p).
$$
Since $F$ is finite and $p,q\in L_0$, this bounds $v$.
It follows from $m=q^{-1}\gamma p\exp(v)$ that $m$ belongs to a fixed compact set.
Increasing $C$ proves all the required estimates for $\gamma\in F$ as well.
\end{proof}

The proof below uses properness arguments similar to those in \cite[Section 9]{KOW} and \cite[Section 7]{BLW}. 
We include a proof for completeness, for an arbitrary $(\Ga,\te)$-proper linear form.
 
\begin{theorem}[Domain of discontinuity]\label{theorem.dd}
The domain $\widetilde{\mathcal M}$ in \eqref{eq.om} has the following properties.
\begin{enumerate}
\item $\widetilde{\mathcal M}$ is open, contains $\widetilde\Omega$, and is invariant
under the left $\Ga$-action and the right $A_\te$-action. The $\Ga$-action on $\widetilde{\mathcal
M}$ is proper.
\item For every $(\Ga,\te)$-proper linear form
$\psi\in\fa_\te^*$, the action $(\gamma, b):x\mapsto \gamma x b^{-1}$ of $\Ga\times \exp(\ker\psi)$
on $\widetilde{\mathcal M}$ is proper, and the induced $\exp(\ker\psi)$-action on $\Ga\ba \widetilde{\mathcal M}$ is
free and proper. In particular, this holds for every hyperplane $\ker\psi$ with $\ker\psi\cap\cal
L_\te=\{0\}$.
\item For each $\psi$ in (2), the right $A_\te$-action induces a complete smooth
$A_\te/\exp(\ker\psi)\simeq\bb R$-flow on the space $\Gamma\backslash \widetilde{\mathcal M}/\exp(\ker\psi)$.
\end{enumerate}
Both quotients $\Gamma\backslash\widetilde{\mathcal M}$ and
$\Gamma\backslash\widetilde{\mathcal M}/\exp(\ker\psi)$
are smooth manifolds.
\end{theorem}
 
\begin{proof} Since $\F_\theta^{(2)}$ is open in $\F_\theta\times\F_{\i\theta}$,
the set of pairs which are not in general position is closed.
Compactness of $\La_\te$ and continuity of $\xi\mapsto \xi^*$ therefore imply that $\widetilde{\mathcal M}$ is open. Its left
$\Ga$-invariance follows from the equivariance of $\xi\mapsto \xi^*$, and its right $A_\te$-invariance
follows from invariance of the two endpoints under $A_\te$.
 
To check $\widetilde\Omega\subset\widetilde{\mathcal M}$, write a pair in $\La_\te^{(2)}$ as $(\xi,\eta^*)$ with $\xi\ne\eta$. For any $\zeta\in\La_\te$, either $\zeta\ne\xi$, in which case $(\xi,\zeta^*)\in\cal F_{\te}^{(2)}$, or $\zeta=\xi$, in which case $(\zeta,\eta^*)\in\cal F_{\te}^{(2)}$. 
This proves the inclusion.
 
Let $L_0\subset G$ be a compact set with $L_0.o\subset\widetilde{\mathcal M}$. 
We claim that $\Ga_0:=\{\ga\in\Ga:\gamma L_0\cap L_0 M_\te\ne\emptyset\}$ is finite. 
Applying Lemma \ref{lem.return} with $v=0$, for all $\ga\in \Ga_0$, either $\mu_\te(\gamma)$ or $\mu_\te(\gamma^{-1})$ is bounded.
Since $\mu_\te|_\Ga$ is proper, the claim follows, proving (1).
 
For part (2), suppose that $\gamma_np_n\exp(-v_n)M_\te=q_n M_\te$ with $p_n,q_n\in L_0$ and $v_n\in\ker\psi$. 
By \eqref{eq.return1}, after passing to a
subsequence, either
$$
 \psi(\mu_\te(\gamma_n))=O(1) \quad\text{or}\quad \psi(\mu_\te(\gamma_n^{-1}))=O(1).
$$
Since $\psi$ is $(\Ga,\te)$-proper, $\gamma_n$ belongs to a finite set. 
Formula \eqref{eq.hopfa} then bounds $v_n$. This proves properness of the $\Ga\times \exp(\ker\psi)$-action and hence of the induced $\exp(\ker\psi)$-action on $\mathcal M$.
The latter action is free because its stabilizers are compact subgroups of $\exp (\ker\psi)\simeq\bb R^{\dim\fa_\te-1}$.
 
If $\ker\psi \cap\cal L_\te=\{0\}$, replacing $\psi$ by $-\psi$ if necessary, we may assume that $\psi>0$ on $\cal L_\te\setminus\{0\}$. 
Such a form $\psi$ is $(\Ga,\te)$-proper \cite[Lemma 4.3]{KOW}, so this case follows from the preceding argument. 
 
Finally, (3) holds since the right $A_\te$-action commutes with $\exp(\ker\psi)$. Since $\Ga$ is torsion-free, both actions are free and proper, so both quotients are smooth manifolds.
\end{proof}
 
We fix the domain $\widetilde{\mathcal M}$ in \eqref{eq.om} and set
$$
 \cal M:=\Ga\ba\widetilde{\mathcal M}, \qquad \Omega:=\Ga\ba\widetilde\Omega.
$$
 
We now apply Lemma \ref{lem.return} to correlations of compactly supported functions.
\begin{proposition}
Let $\mathsf m$ be a locally finite Borel measure on $\mathcal M$.
For any $u\in\mathfrak a_\theta^+\setminus\mathcal L_\theta$ and $f_1,f_2\in C_c(\mathcal M)$, there exists $T_0>0$ such that
$$
 \int_{\mathcal M} f_1(x\exp(tu))f_2(x)\,d\mathsf m(x)=0 \qquad\text{for all }t\geq T_0.
$$
\end{proposition}

\begin{proof} Fix $u\in\mathfrak a_\theta^+\setminus\mathcal L_\theta$ and let $\cal C\subset\fa_\te$ be
a closed cone with $\cal L_\te\setminus\{0\}\subset\op{int}\cal C$ and $u\notin \cal C$.
We claim that for every compact subset $L_0\subset G$ with $L_0.o\subset\widetilde{\cal M}$,
$\Ga  L_0\exp(v)\cap L_0M_\te=\emptyset$ for all sufficiently large $v\in\fa_\te^+\setminus\cal C$. 
Choose $L_0$ whose image in $\mathcal M$ contains $\supp f_1\cup\supp f_2$. The proposition then follows from the claim.
Suppose the claim does not hold.
Then there are $p_n,q_n\in L_0$, $m_n\in M_\te$, $\gamma_n\in\Ga$, and $v_n\in\fa_\te^+\setminus\cal C$ such that
$$
 t_n:=\|v_n\|\to\infty, \qquad \ga_n p_n\exp(v_n)=q_n m_n.
$$
Applying Lemma \ref{lem.return} to $\gamma_n p_n\exp(v_n)=q_nm_n$ gives
$\mu_\te(\gamma_n^{-1})=v_n+O_{L_0}(1)$. 
Pass to a subsequence with $\frac{v_n}{t_n}\to v_\infty$. 
Then $\|v_\infty\|=1$ and $\frac{\mu_\te(\gamma_n^{-1})}{t_n}\to v_\infty$. 
Definition \ref{def.lc} implies $v_\infty\in\cal L_\te$. 
But since $\frac{v_n}{t_n}\notin\cal C$, we have $v_\infty\notin\op{int}\cal C$, which is a contradiction.
\end{proof}

\section{Product formulas for the BMS, BR and Haar measures}\label{sec.Norb}
\noindent \textbf{The generalized BMS measure: $m_{\nu_1,\nu_2}$.}
Fix linear forms $\psi_1\in \fa_\te^*$, $\psi_2\in \mathfrak a_{\i\te}^*$. 
Let $\nu_1=\nu_{\psi_1}$ and $\nu_2=\nu_{\psi_2}$ be, respectively, $(\Gamma,\psi_1)$ and $(\Gamma,\psi_2)$ conformal measures on $\F_\te$ and $\F_{\i\te}$. 
Using the Hopf parametrization, define a locally finite Borel measure $\tilde m_{\nu_1, \nu_2}$ on $G/M_\te$ as follows: for $gM_\te=(g^+, g^-, b)\in \F_\te^{(2)}\times \mathfrak a_\te$,
\begin{equation}\label{eq.BMS0}
 d\tilde m_{\nu_1, \nu_2} (gM_\te)=e^{\psi_1 (\beta_{g^+}^\te(e, g))+\psi_2( \beta_{g^-}^{\i\te} (e, g )) } \; d\nu_{1} (g^+) d\nu_{2}(g^-) db,
\end{equation}
where $db=d\ell(b)$ is the Lebesgue measure on $\mathfrak a_\te$.
\begin{lemma}\label{qa}
The measure $\tilde m_{\nu_1, \nu_2}$ is left $\Gamma$-invariant and right $A_\te$-quasi-invariant: for every $a\in A_\te$,
$$
 a_*\tilde m_{\nu_1, \nu_2}=e^{(-\psi_1+\psi_2\circ \i)(\log a)}\,\tilde m_{\nu_1, \nu_2}
$$
where $a_*$ denotes push-forward under $gM_\te\mapsto gaM_\te$.
\end{lemma}
\begin{proof}
For $\gamma\in\Ga$ and $\xi\in\F_\te$,  \eqref{eq.basic0} and Definition \ref{def.conf} give
$$
 \beta_{\gamma\xi}^\te(e,\gamma g) =\beta_{\gamma\xi}^\te(e,\gamma)+\beta_\xi^\te(e,g),\qquad d\nu_1(\gamma\xi) =e^{-\psi_1(\beta_{\gamma\xi}^\te(e,\gamma))}\,d\nu_1(\xi).
$$
The same identities hold for $\nu_2$ with $\i\te$ and $\psi_2$.
The exponential factors cancel in \eqref{eq.BMS0}, while the Hopf coordinate changes by a translation. 
This proves left $\Ga$-invariance, as in \cite[Lemma 3.6]{ELO}.
For $a\in A_\te$, \eqref{eq.basic1} gives
$$
 \beta_{g^+}^\te(e,ga)=\beta_{g^+}^\te(e,g)+\log a,\qquad \beta_{g^-}^{\i\te}(e,ga) =\beta_{g^-}^{\i\te}(e,g)-\i\log a.
$$
The endpoints are unchanged and $b$ is translated by $\log a$. Thus
$$
 d\tilde m_{\nu_1,\nu_2}(gaM_\te) =e^{(\psi_1-\psi_2\circ\i)(\log a)} d\tilde m_{\nu_1,\nu_2}(gM_\te),
$$
which gives the stated push-forward formula.
\end{proof}
 
Restrict $\tilde m_{\nu_1,\nu_2}$ to $\widetilde{\mathcal M}$ and denote the induced quotient measure on $\mathcal M$ by $m_{\nu_1,\nu_2}$. 
We call $m_{\nu_1,\nu_2}$ the generalized BMS measure associated to $(\nu_1,\nu_2)$.
 
\medskip
\noindent
\textbf{BMS measure: $m_{\nu_\psi,\nu_{\psi\circ\i}}^{\BMS}$.} Let $\psi\in\fa_\te^*$ and let $\nu_\psi$ and $\nu_{\psi\circ\i}$ be, respectively, $(\Gamma,\psi)$- and $(\Gamma,\psi\circ\i)$-PS measures on $\F_\te$ and $\F_{\i\te}$. 
We set
$$
 m^{\BMS}_{\nu_\psi,\nu_{\psi\circ\i}} :=m_{\nu_\psi,\nu_{\psi\circ\i}}
$$
and call it the Bowen-Margulis-Sullivan measure associated to $(\nu_\psi,\nu_{\psi\circ\i})$. 
By Lemma \ref{qa}, this measure is $A_\te$-invariant. Writing $[g]=\Ga gM_\te\in\mathcal M$, its support is
$
 \operatorname{supp}(m^{\BMS}_{\nu_\psi,\nu_{\psi\circ\i}}) =\{[g]\in\mathcal M:g^+\in\La_\te,\ g^-\in\La_{\i\te}\}.
$
 
\medskip
\noindent
\textbf{Burger-Roblin measures.} We define the BR and dual BR measures, respectively, by
$$
 m^{\BR}_{\nu_{\psi\circ\i}} :=m_{m_\te,\nu_{\psi\circ\i}},\qquad m^{\BR_*}_{\nu_\psi}:=m_{\nu_\psi,m_{\i\te}}.
$$
 
These generalize the rank-one measures of Burger \cite{Bur} and Roblin \cite{Ro}.
Their supports are
\begin{align*}
 \operatorname{supp}(m^{\BR}_{\nu_{\psi\circ\i}})&=\{[g]\in\mathcal M:g^-\in\La_{\i\te}\},\quad 
 \operatorname{supp}(m^{\BR_*}_{\nu_\psi})&=\{[g]\in\mathcal M:g^+\in\La_\te\}.
\end{align*}
 
\medskip
\noindent
\textbf{Haar measure $m^{\Haar}$.} Recall that the $K$-invariant probability measure $m_\te$ is a 
$(G, 2\rho_\te)$-conformal measure.
We denote by $dx=dm^{\Haar}(x)$ the generalized BMS measure associated to the pair $(m_\te,m_{\i\te})$:
$$
 dx=dm^{\Haar}:=dm_{m_\te,m_{\i\te}}.
$$
The corresponding measure on $G/M_\te$ is left $G$-invariant: apply the proof of Lemma \ref{qa} with $G$ in place of $\Ga$.
 
For the rest of the section,  fix a linear form $\psi\in \fa_\theta^*$, a $(\Gamma, \psi)$-PS measure $\nu_\psi$ on $\cal F_\te$ and a $(\Gamma, \psi\circ \i)$-PS measure $\nu_{\psi\circ \i}$ on $\F_{\i\te}$. 
We write
$$
 \nu:=\nu_{\psi},\;\nu_{\op{i}}:=\nu_{\psi\circ \op{i}}, \; \tilde m^{\BMS}:=\tilde m^{\BMS}_{\nu_\psi, \nu_{\psi\circ \i}}, \; \tilde m^{\BR}:=\tilde m^{\BR}_{\nu_{\psi\circ \i}}, \; \tilde m^{\BR_*}:=\tilde m^{\BR_*}_{\nu_\psi}.
$$
\subsection*{PS measures on $N_\te^\pm$}
For $g\in G$, define locally finite measures on $N_\te^+$ and $N_\te^-$ by
\begin{align*}
 d\mu_{g,N_\te^+}^{\PS}(n)&=e^{\psi(\beta_{(gn)^+}^\te(e,gn))}\,d\nu((gn)^+),\\
 d\mu_{g,N_\te^-}^{\PS}(h)&=e^{\psi\circ\i(\beta_{(gh)^-}^{\i\te}(e,gh))}\,d\nu_{\op{i}}((gh)^-).
\end{align*}
These measures are locally finite, since $n\mapsto(gn)^+$ and $h\mapsto(gh)^-$ are diffeomorphisms onto open subsets of $\F_\te$ and $\F_{\i\te}$, the measures $\nu$ and $\nu_{\op{i}}$ are finite, and the densities are continuous.
 
Recalling that $A_\te M_\te$ normalizes $N_\te^{\pm}$, we will use the following lemma:
 
\begin{lemma}\label{MANconj} For any $g\in G,\,a\in A_\te,\,n_0,n\in N_\te^+$ and $m\in M_\te$, we have
\begin{equation*}
 d\mu_{g,N_\te^+}^{\mathrm{PS}}(amn_0n m^{-1}a^{-1})=e^{-\psi(\log a)}d\mu_{gam,N_\te^+}^{\mathrm{PS}}(n_0n).
\end{equation*}
\end{lemma}
\begin{proof} 
By \eqref{eq.basic0} and \eqref{eq.basic1}, we have $(gamn_0nm^{-1}a^{-1})^+=gamn_0n^+$, and
\begin{align*}
 &\beta_{gamn_0n^+}^{\te}(e,gamn_0nm^{-1}a^{-1})=\beta_{gamn_0n^+}^\te(e,gamn_0n)+\beta_{e^+}^\te(e,m^{-1}a^{-1}).
\end{align*}
Also note that $\beta_{e^+}^\te(e,m^{-1}a^{-1})=-\log a$ by \eqref{eq.basic1}.
Hence
\begin{align*}
 d\mu_{g,N_\te^+}^{\mathrm{PS}}(amn_0nm^{-1} a^{-1})
 &=e^{-\psi(\log a)}e^{\psi(\beta_{gamn_0n^+}^\te(e,gamn_0n))}d\nu (gamn_0n^+) \\
 &=e^{-\psi(\log a)}d\mu_{gam,N_\te^+}^{\mathrm{PS}}(n_0n).
\end{align*}
\end{proof}
 
\subsection*{Product formulas for BMS measures}

The product map $N_\te^+\times P_\te^-\to G$ is a diffeomorphism onto an open neighborhood of $e$. Hence we have the following product formula.
\begin{lemma}\label{NPdecomp}
For $g\in G$, $f\in C_c(gN_\te^+ N_\te^-A_\te.o)$,
\begin{equation*}
 \tilde{m}^{\mathrm{BMS}}(f)=\int_{N_\te^+}\left(\int_{N_\te^-}\int_{A_\theta} f(gnha.o)\, da \; d\mu_{gn,N_\te^-}^{\mathrm{PS}}(h)\right)d\mu_{g,N_\te^+}^{\mathrm{PS}}(n).
\end{equation*}
\end{lemma}
\begin{proof}
By \eqref{eq.basic0} and \eqref{eq.basic1}, $(gnha)^-=(gnh)^-$, $(gnha)^+=(gn)^+$, and
$$
\begin{aligned}
 \beta_{(gn)^+}^{\te}(e,gnha)&=\beta_{(gn)^+}^{\te}(e,gn)+\log a,\\
 \beta_{(gnh)^-}^{\i\te}(e,gnha)&=\beta_{(gnh)^-}^{\i\te}(e,gnh)-\i\log a.
\end{aligned}
$$
Hence
\begin{align*}
 d\tilde{m}^{\mathrm{BMS}}(gnha.o)
 =&e^{\psi( \op{i}(\beta_{gnh^-}^{\i\te}(e,gnh))+\beta_{gnh^+}^{\te}(e,gnh))} d\nu_{\op{i}}(gnh^-)\, d\nu(gn^+)\,da\\
 =&da\,d\mu_{gn,N_\te^-}^{\mathrm{PS}}(h)\,d\mu_{g,N_\te^+}^{\mathrm{PS}}(n),
\end{align*}
proving the claim.
\end{proof}
 
Similarly, using $N_\te^-A_\te N_\te^+$-coordinates, we have:
\begin{lemma}\label{PNdecmop}
For $g\in G$, $f\in C_c(gN_\te^-A_\te N_\te^+.o)$, 
\begin{equation*}
 \tilde{m}^{\mathrm{BMS}}(f)=\int_{N_\te^-}\int_{A_\te}\left(\int_{N_\te^+}f(ghan.o)\,d\mu_{gha,N_\te^+}^{\mathrm{PS}}(n)\right)e^{-\psi(\log a)}\,da\,d\mu_{g,N_\te^-}^{\mathrm{PS}}(h).
\end{equation*}
\end{lemma}
 
Define, for $ha\in N_\te^-A_\te$,
$$
 d\mu_{g,N_\te^-A_\te}^{\mathrm{PS}}(ha)=e^{-\psi(\log a)}\,da\,d\mu_{g,N_\te^-}^{\mathrm{PS}}(h).
$$
The formula in Lemma \ref{PNdecmop} becomes: for any $ f\in C_c(gN_\te^-A_\te N_\te^+.o)$,
\begin{equation}\label{PNgood}
 \tilde{m}^{\mathrm{BMS}}(f)=\int_{N_\te^-A_\te}\int_{N_\te^+} f(gha n.o)\,d\mu_{gha,N_\te^+}^{\mathrm{PS}}(n)\,d\mu_{g,N_\te^-A_\te}^{\mathrm{PS}}(ha).
\end{equation}
 
\noindent\textbf{Lebesgue measures on $gN_\te^{\pm}.o$.} For $g\in G$,
we note that the Haar measure on $N_\te^{\pm}$ can be given as follows: for $n\in N_\te^{\pm}$,
$$
 d\mu_{g,N_\te^-}^{\mathrm{Leb}}(n)= e^{2\rho_{\i\te} (\beta_{gn^-}^{\i\te}(e, gn)) }dm_{\i\te}(gn^-);
$$
\be\label{eq.Leb2}
 d\mu_{g,N_\te^+}^{\mathrm{Leb}}(n)=e^{2\rho_\te (\beta_{gn^+}^\te (e,gn) ) }dm_{\te} (gn^+).
\ee
The $(G,2\rho_\te)$-conformality of $m_\te$ \eqref{gc0}, the identity \eqref{eq.basic0}, and \eqref{eq.Leb2} give
$$
 d\mu_{g,N_\te^+}^{\mathrm{Leb}}(n) =e^{2\rho_\te(\beta_{n^+}^{\te}(e,n))}\,dm_\te(n^+) =d\mu_{e,N_\te^+}^{\mathrm{Leb}}(n).
$$
It follows that this measure is left $N_\te^+$-invariant. 
The same argument applies to $N_\te^-$. We write these Haar measures as $dn$. 
We write $da=d\ell(\log a)$ for the Haar measure on $A_\theta$ corresponding to the Lebesgue measure in \eqref{eq.BMS0}.
 
\subsection*{Product formulas for BR measures}
The following formulas can be verified as in the proof of Lemma \ref{PNdecmop}: for all $f\in C_c(gN_\te^-A_\te N_\te^+.o)$,
\begin{equation}\label{BRPN}
 \tilde m^{\BR}(f) =\int_{N_\te^-A_\te}\int_{N_\te^+} f(ghan.o)\,dn\,d\mu_{g,N_\te^-A_\te}^{\PS}(ha).
\end{equation}
For all $f\in C_c(gN_\te^+A_\te N_\te^-.o)$,
$$
 \tilde{m}^{\mathrm{BR}_*}(f)=\int_{N_\te^+}\int_{A_\te}\int_{N_\te^-}f(gnah.o)\,dh\,e^{\psi(\log a)}\,da\,d\mu_{g,N_\te^+}^{\mathrm{PS}}(n).
$$
For all $f\in C_c(gN_\te^-A_\te N_\te^+.o)$, we also have
\begin{equation}\label{BRPNdual}
 \tilde{m}^{\mathrm{BR}_*}(f)=\int_{N_\te^-}\int_{A_\te}\int_{N_\te^+}f(ghan.o)\,d\mu_{gha,N_\te^+}^{\mathrm{PS}}(n)e^{-2\rho_\te(\log a)}\,da\,dh.
\end{equation}
We also have the following description of the BR measures on $G/M_\te$:
\begin{lemma}\label{lem.BRD}
Let $d\ell$ denote the Lebesgue measure on $\mathfrak a_\te$.
We have:
\begin{itemize}
\item
For $g=k\exp(H)nm\in K A_\te N_\te^+M_\te$,
\begin{equation*}
 d\tilde m_{\nu_{\psi\circ\i}}^{\BR}(gM_\te)=e^{-\psi (H)}\,dn\,d\ell(H)\,d\nu_{\psi\circ\i} (k^-).
\end{equation*}
\item
For $g=k\exp(H)nm\in KA_\te N_\te^-M_\te$,
\begin{equation*}
 d\tilde m_{\nu_\psi }^{\BR_*}(gM_\te)=e^{\psi(H)}\,dn\,d\ell(H)\,d\nu_{\psi} (k^+).
\end{equation*}
\end{itemize}
\end{lemma}
\begin{proof}
For $g=k\exp(H)nm\in K A_\te N_\te^+ M_\te$, we have $\beta_{g^-}^{\i\te}(e,g)=-\i (H)$ by
\eqref{eq.basic1}. 
Since $m_\te$ is a $(G, 2\rho_\te)$-conformal measure, we have
\begin{align*}
 &d\tilde m^{\BR}_{\nu_{\psi\circ\i}} (gM_\te) \\ &=e^{2\rho_\te(\beta_{k\exp(H)n^+}^\te(e,k\exp (H)n)) }e^{-\psi \circ\i(\i(H))} \, dm_\te(k\exp(H)n^+)\,d\ell(H) \, d \nu_{\psi\circ\i}(k^-)\\
 &=e^{-\psi(H)}\,dn\,d\ell(H)\, d \nu_{\psi\circ\i}(k^-).
\end{align*}
The second formula follows similarly.
\end{proof}

\subsection*{Product formulas for Haar measure}
Recall that $\tilde m^{\Haar}=\tilde m_{m_\te,m_{\i\te}}$ is the $G$-invariant measure on $G/M_\te$. 
The same computation as in the proof of Lemma \ref{lem.BRD}, together with $2\rho_{\i\te}\circ\i=2\rho_\te$, gives
\begin{align*}
 d\tilde m^{\Haar}(k\exp(H)nM_\te)&=e^{2\rho_\te(H)}\,dn\,d\ell(H)\,dm_\te(k^+), && n\in N_\te^-,\\
 d\tilde m^{\Haar}(k\exp(H)nM_\te)&=e^{-2\rho_\te(H)}\,dn\,d\ell(H)\,dm_{\i\te}(k^-), && n\in N_\te^+.
\end{align*} 
Similarly, by \eqref{BRPNdual}, for every $f\in C_c(gN_\te^-A_\te N_\te^+.o)$,
$$
 \tilde m^{\Haar}(f)=\int_{N_\te^-}\int_{A_\te}\int_{N_\te^+}f(ghan.o)e^{-2\rho_\te(\log a)}\,dn\,da\,dh.
$$
\section{Local coordinates and smoothing}\label{sec.lcd}
Recall the distance $d$ on $G$ fixed in Section \ref{ps}. Choose a smooth Riemannian metric on $\mathcal M$, with distance $d_{\mathcal M}$, and equip $\widetilde{\mathcal M}$ with its pullback.

We fix $C^k$-norms on $C_c^k(\mathcal M)$ and $C_c^k(\widetilde{\mathcal M})$ defined using derivatives of order at most $k$ in smooth coordinates. On functions supported in a fixed compact set, any two such choices are equivalent.
 
For $\epsilon>0$, let $G_\epsilon:=\{g\in G: d(e,g)<\epsilon\}$.
For a closed subgroup $S<G$, set
$$
 S_\epsilon:=S\cap G_\epsilon.
$$
Let $\norm{\cdot}_{C^k(S)}$ be the $C^k$-norm on $S$ defined using a fixed basis of left-invariant vector fields.
The following construction is standard; see \cite[Section 2, p. 138]{Friedrichs} for the Euclidean case.
\begin{lemma}\label{lem.bump}\rm
For any closed subgroup $S<G$ and $0<\epsilon<1$, there exists a nonnegative $\eta_\epsilon\in C_c^\infty(S_\epsilon)$ such that
$$
 \int_S\eta_\epsilon(h)\,dh=1,\qquad
 \|\eta_\epsilon\|_{C^k(S)}
 \ll_k\epsilon^{-(k+\op{dim}S)}
 \quad\text{for all} \quad k\geq0.
$$
\end{lemma}
 
By \cite[Proposition 8.45]{Knapp}, the set
$$
 \mathcal O:=N_\theta^+N_\theta^-A_\theta.o
$$
is open and dense in $G/M_\theta$, and, for every $g\in G$, the map
$$
 N_\theta^+\times N_\theta^-\times A_\theta\to g\mathcal O,\qquad (n,h,a)\mapsto gnha.o
$$
is a diffeomorphism. Using exponential coordinates on $N_\theta^+$, $N_\theta^-$, $A_\theta$, we define the $C^k$-norms $\|\cdot\|_{C^k(g\mathcal O)}$. 
On functions supported in a fixed compact subset of $g\mathcal O\cap\widetilde{\mathcal M}$, these norms are equivalent to the chosen norm $\|\cdot\|_{C^k(\widetilde{\mathcal M})}$.
 
\subsection*{Smoothing of continuous functions}
For $f\in C_c(\mathcal M)$ and $\epsilon>0$, define
$$
 f_\epsilon^+(x):=\sup\{f(y):d_{\mathcal M}(x,y)\leq\epsilon\},\qquad f_\epsilon^-(x):=\inf\{f(y):d_{\mathcal M}(x,y)\leq\epsilon\}.
$$
Apply Lemma \ref{lem.bump} with $S=G$ to choose $\eta_\epsilon\in C_c^\infty(G_\epsilon)$ such that
\begin{equation}\label{eq.eta}
 \int_G \eta_\epsilon(h)\,dh=1,\quad\text{ and }\quad \norm{\eta_\epsilon}_{C^k(G)}\ll_k \epsilon^{-(k+\op{dim}G)}
\end{equation}
for all $k\geq 0$. Choose relatively compact coordinate neighborhoods $U_j\subset\mathcal M$ and a
locally finite smooth partition of unity $\{\chi_j\}$ with $\supp\chi_j\subset U_j$, and smooth
local sections $s_j:U_j\to G$ satisfying $[s_j(x)]=x$. Choose $0<c_j<1$ so small that
$$
 d_{\mathcal M}(x,[s_j(x)h])\leq\epsilon \quad\text{for }x\in\supp\chi_j,\ h\in G_{c_j\epsilon}, \quad 0<\epsilon<1.
$$
In the following integrals, functions of $[q]$ are extended by zero when
$qM_\te\notin\widetilde{\mathcal M}$. 
Define
\begin{equation}\label{eq.ffe}
    F_{\epsilon}^\pm(x) :=\sum_j\chi_j(x)\int_G f_{2\epsilon}^\pm([s_j(x)h]) \eta_{c_j\epsilon}(h^{-1})\,dh,
\end{equation}
where each summand is extended by zero outside $U_j$. The supports of these functions lie in
$
 \{x:d_{\mathcal M}(x,\supp f)\leq3\epsilon\}.
$
For sufficiently small $\epsilon>0$, this set is contained in a fixed relatively compact neighborhood of $\supp f$. 
The functions are smooth: in the $j$th integral, the substitution $z=s_j(x)h$ gives
$$
 \int_G f_{2\epsilon}^\pm([z]) \eta_{c_j\epsilon}(z^{-1}s_j(x))\,dz,
$$
so derivatives in $x$ act on the smooth factor $\eta_{c_j\epsilon}(z^{-1}s_j(x))$. 
 
\begin{lemma}\label{Lem.FE1}
For all sufficiently small $\epsilon>0$ and all $x,y\in\mathcal M$ with $d_{\mathcal
M}(x,y)\leq\epsilon$,
$$
 F_{\epsilon}^-(x)\leq f(y)\leq F_{\epsilon}^+(x).
$$
\end{lemma}
 
\begin{proof}
If $\chi_j(x)\ne0$ and $h\in G_{c_j\epsilon}$, set $z=[s_j(x)h]$. Then
$$
 d_{\mathcal M}(z,y) \leq d_{\mathcal M}(z,x)+d_{\mathcal M}(x,y) \leq2\epsilon,
$$
and hence $ f_{2\epsilon}^-(z)\leq f(y)\leq f_{2\epsilon}^+(z)$.
Integrate against $\eta_{c_j\epsilon}(h^{-1})\,dh$, multiply by $\chi_j(x)$, and sum over $j$. The conclusion follows from $\sum_j\chi_j(x)=1$ and $\int_G\eta_{c_j\epsilon}(h^{-1})\,dh=1$, since $G$ is unimodular.
\end{proof}
 
\begin{lemma}\label{lem.e1}
Let $\cal K\subset\mathcal M$ be a compact set. For all sufficiently small $\epsilon>0$ depending on $\cal K$, all $k\in\mathbb{N}$ and real-valued $f\in
C_c^\infty (\cal K)$,
\begin{enumerate}
\item
$\norm{F_{\epsilon}^\pm}_{C^k(\mathcal M)}\ll_{\cal K, k}
\epsilon^{-k}\norm{f}_{\infty}.$
\item $m^{\mathrm{BMS}}(F_{\epsilon}^+-F_{\epsilon}^-)\ll_{\cal K} \epsilon \norm{f}_{C^1(\cal M)}$.
\end{enumerate}
\end{lemma}
\begin{proof}
By \eqref{eq.eta}, the integral of the absolute value of each derivative of order at most $k$ of $\eta_{c_j\epsilon}$ is $O_{j,k}(\epsilon^{-k})$: its support has volume $O_j(\epsilon^{\dim G})$.
Differentiate the preceding expression with $\eta_{c_j\epsilon}(z^{-1}s_j(x))$.
The derivatives of the finitely many sections $s_j$ and functions $\chi_j$ involved are bounded on the fixed compact neighborhood.
Since $|f_{2\epsilon}^\pm|\leq\norm{f}_\infty$, this gives (1).
The mean value estimate on the fixed compact neighborhood gives $\norm{f_{2\epsilon}^+-f_{2\epsilon}^-}_\infty\ll_{\cal K} \epsilon \norm{f}_{C^1(\cal M)}$, and hence
$
 \|F_{\epsilon}^+-F_{\epsilon}^-\|_{\infty} \ll_{\cal K}\epsilon\|f\|_{C^1(\mathcal M)}.
$
Since $m^{\BMS}$ is locally finite and the supports lie in a fixed compact set, (2) follows.
\end{proof}
 
\section{BMS local mixing and translates of PS measures}\label{sec.set}
We fix a vector $u\in\mathcal L_\theta\cap\operatorname{int}\mathfrak a_\theta^+$ and a $(\Gamma,\theta)$-proper tangent form $\psi\in\mathfrak a_\te^*$  with $\psi(u)=\psi_\Gamma^\theta(u)>0$. 
Let $\nu=\nu_\psi$ and $\nu_{\i}=\nu_{\psi\circ\i}$ be $(\Gamma,\psi)$- and $(\Gamma,\psi\circ\i)$-PS measures on $\F_\te$ and $\F_{\i\te}$, respectively.
As before, we set
$$
 m^{\BMS}= m^{\BMS}_{\nu, \nu_{\i}}, \quad \mu^{\PS}_{g,N_\te^+}=\mu^{\PS}_{g,N_\te^+, \nu},\quad \mu^{\PS}_{g,N_\te^-}=\mu^{\PS}_{g,N_\te^-, \nu_{\i}}.
$$
For $t>0$ and $v\in\ker\psi$, we use $a_{tu+\sqrt t\,v}=\exp(tu+\sqrt t\,v)\in A_\te$.
Fix a function $\Psi:(0,\infty)\to(0,\infty)$ and a bounded Borel function $J:\ker\psi\to(0,\infty)$.
\begin{Def}\label{lmp}\rm
\begin{enumerate}
\item We say that $m^{\BMS}$ satisfies \emph{local mixing} with the pair $(\Psi,J)$ if, for every $v\in\ker\psi$ and $f_1,f_2\in C_c(\mathcal M)$,
$$
\lim_{t\to \infty} \Psi(t)\int_{\mathcal M}f_1(xa_{tu+\sqrt t\,v})f_2(x)\,dm^{\BMS}(x)=J(v)m^{\BMS}(f_1)m^{\BMS}(f_2).
$$
 
\item We say that the \emph{uniform bound} holds if, for every $f_1,f_2\in C_c(\mathcal M)$, there exists $C=C(f_1,f_2)>0$ such that, for all $t\geq1$ and $v\in\ker\psi$ with $a_{tu+\sqrt t\,v}\in A_\te^+$,
$$
 \left|\Psi(t)\int_{\mathcal M}f_1(xa_{tu+\sqrt t\,v})f_2(x)\,dm^{\BMS}(x)\right|\leq C.
$$
\end{enumerate}
\end{Def}
Set
\begin{equation*}
\begin{split}
 \err(a_{tu+\sqrt t\,v},f_1,f_2)&:=\Psi(t)\int_{\mathcal M}f_1(xa_{tu+\sqrt t\,v})f_2(x)\,dm^{\BMS}(x)\\
 &\quad{}-J(v)m^{\BMS}(f_1)m^{\BMS}(f_2).
\end{split}
\end{equation*}
 
\begin{Def}[Polynomial rates]\label{def.em}\rm
For $d>0$, consider the following conditions, each with some integer $k\geq1$. 
\begin{enumerate}
\item For every compact $\mathcal K\subset\mathcal M$ and fixed $v\in\ker\psi$,
\begin{equation}\label{eq.kd}
 \err(a_{tu+\sqrt t\,v},f_1,f_2)=O_{\mathcal K,v}\!\left(t^{-d}\|f_1\|_{C^k(\mathcal M)}\|f_2\|_{C^k(\mathcal M)}\right)
\end{equation}
as $t\to\infty$, uniformly for $f_1,f_2\in C_c^\infty(\mathcal K)$.
\item We have $J\in L^1(\ker\psi)$, and for every compact $\mathcal K\subset\mathcal M$ there is a nonnegative bounded Borel function $Q_{\mathcal K}\in L^1(\ker\psi)$ such that
\begin{equation}\label{eq.kd2}
 |\err(a_{tu+\sqrt t\,v},f_1,f_2)|\leq Q_{\mathcal K}(v)t^{-d}\|f_1\|_{C^k(\mathcal M)}\|f_2\|_{C^k(\mathcal M)}
\end{equation}
for all $f_1,f_2\in C_c^\infty(\mathcal K)$, $t\geq1$ and $v\in\ker\psi$ with $a_{tu+\sqrt t\,v}\in A_\theta^+$.
\end{enumerate}
\end{Def}
 
The normalizing function in local mixing has the following properties.
\begin{lemma}\label{lem.ratio}\rm
Suppose that $m^{\BMS}$ satisfies local mixing as in Definition \ref{lmp}.
\begin{enumerate}
\item For every $C>0$,
$$
 \sup_{|s|\leq C}\left|\frac{\Psi(T+s)}{\Psi(T)}-1\right|\to0\quad\text{ as }\quad T\to\infty.
$$
\item For every $\eta>0$, there exists $t_\eta>0$ such that
$$
 e^{-\eta(T-t+1)}\leq\frac{\Psi(T)}{\Psi(t)}\leq e^{\eta(T-t+1)}\quad\text{ for all }\quad t_\eta\leq t\leq T.
$$
In particular, $$\log\Psi(T)=o(T).$$
\item Suppose, in addition, that the polynomial estimate \eqref{eq.kd} in Definition \ref{def.em}(1) holds with order $d>0$. Then there exist $C_\Psi,T_0>0$ such that
$$
 \left|\log\frac{\Psi(T+s)}{\Psi(T)}\right|\leq C_\Psi(|s|+1)T^{-d}\quad\text{ for all }\quad T\geq T_0,\ |s|\leq T/2.
$$
\end{enumerate}
\end{lemma}
 
\begin{proof}
Fix a nonnegative $f\in C_c^\infty(\mathcal M)$ with $m^{\BMS}(f)>0$, and set
$$
 f_s(x):=f(xa_{su}),\qquad L_t(\phi):=\Psi(t)\int_{\mathcal M}\phi(xa_{tu})f(x)\,dm^{\BMS}(x).
$$
By $A_\theta$-invariance, $m^{\BMS}(f_s)=m^{\BMS}(f)$, and
\begin{equation}\label{eq.ratio}
 L_t(f_s)=\frac{\Psi(t)}{\Psi(t+s)}L_{t+s}(f).
\end{equation}
We prove (1).
Let $\cal K$ be a compact set such that $\supp f_s\subset\cal K$ for all $|s|\leq C$.
Let $\chi\in C_c(\mathcal M)$ be a nonnegative function such that $\chi|_{\mathcal K}=1$. 
By Definition \ref{lmp}, $L_t(\chi)$ is bounded for all sufficiently large $t$. 
Since $f\geq 0$, $|L_t(\phi)|\leq\|\phi\|_\infty L_t(\chi)$ for all $\phi$ supported in $\mathcal K$.
 
By compactness, for every $\epsilon>0$, there exist $s_1,\ldots,s_N\in[-C,C]$ such that, for each $|s|\leq C$, there is an $i$ with $\|f_s-f_{s_i}\|_\infty<\epsilon$. Then
$$
 |L_t(f_s)-L_t(f_{s_i})|\leq\epsilon L_t(\chi).
$$
By Definition \ref{lmp}, $L_t(f_{s_i})\to J(0)m^{\BMS}(f)^2$ for each $1\le i\le N$.  
Therefore $L_t(f_s)\to J(0)m^{\BMS}(f)^2$ uniformly for $|s|\leq C$.  
It follows that
$
 L_{t+s}(f)\to J(0)m^{\BMS}(f)^2>0
$
uniformly for $|s|\leq C$. Dividing in \eqref{eq.ratio} proves (1).

Applying (1) with $C=1$, for every $\eta>0$ there exists $t_\eta>0$ such that
$$
 |\log\Psi(t+h)-\log\Psi(t)|\leq\eta\quad\text{ for all }\quad t\geq t_\eta,\ 0\leq h\leq1.
$$
Splitting $[t,T]$ into at most $T-t+1$ intervals of length at most one and adding these bounds proves (2).
 
To prove (3), we apply Definition \ref{def.em}(1) with $f_h$ at $v=0$.
We obtain
$$
 L_t(f_h)=J(0)m^{\BMS}(f)^2+O(t^{-d}),\qquad L_{t+h}(f)=J(0)m^{\BMS}(f)^2+O(t^{-d})
$$
uniformly in $|h|\leq 1$. 
Hence, \eqref{eq.ratio} implies
$
 |\log\Psi(t+h)-\log\Psi(t)|\ll t^{-d}
$ for all $|h|\le 1$ and 
for all sufficiently large $t$. Splitting the interval between $T$ and $T+s$ into at most $|s|+1$ intervals of length at most one, all contained in $[T/2,3T/2]$, gives
$
 \left|\log\frac{\Psi(T+s)}{\Psi(T)}\right|\ll(|s|+1)T^{-d}.
$
\end{proof}
 
The main result of this section is the following.
\begin{proposition}\label{p1}\rm
Let $g\in G$, $f\in C_c^\infty(\mathcal M)$, and $\phi\in C_c^\infty(N_\theta^+)$ with $g(\supp\phi).o\subset\widetilde{\mathcal M}$. Then (1) and (2) below hold under conditions (1) and (2) of Definition \ref{lmp}, respectively.
\begin{enumerate}
\item For every  $v\in\ker\psi$,
\begin{align*}
 &\Psi(t)\int_{N_\theta^+}f([gn]a_{tu+\sqrt t\,v})\phi(n)\,d\mu_{g,N_\theta^+}^{\PS}(n)\\
 &=J(v)m^{\BMS}(f)\mu_{g,N_\theta^+}^{\PS}(\phi)+\err'(a_{tu+\sqrt t\,v},f,\phi),
\end{align*}
where $\err'(a_{tu+\sqrt t\,v},f,\phi)\to0$ as $t\to\infty$.
\item There exists $C'=C'(g,f,\phi)>0$ such that, for all $t\geq1$ and $v\in\ker\psi$ with $a_{tu+\sqrt t\,v}\in A_\theta^+$,
$$
 \left|\Psi(t)\int_{N_\theta^+}f([gn]a_{tu+\sqrt t\,v})\phi(n)\,d\mu_{g,N_\theta^+}^{\PS}(n)\right|\leq C'.
$$
\end{enumerate}
Assume now that $\Gamma$ is linearly $\theta$-regular \eqref{eq.rc}.
Then (3) and (4) below hold under conditions (1) and (2) of Definition \ref{def.em}, respectively, for some $d'>0$ depending only on $d$, $k$ and $G$.
 
Fix compact sets $\mathcal K_1\subset\mathcal M$ and $\mathcal K_2\subset N_\theta^+$ with $g\mathcal K_2.o\subset\widetilde{\mathcal M}$. Both estimates below are uniform for $f\in C_c^\infty(\mathcal K_1)$ and $\phi\in C_c^\infty(\mathcal K_2)$.
\begin{enumerate}
\item[(3)] For every $v\in\ker\psi$,
$$
 \err'(a_{tu+\sqrt t\,v},f,\phi)=O_{g,\mathcal K_1,\mathcal K_2,v}\!\left(t^{-d'}\|f\|_{C^k(\mathcal M)}\|\phi\|_{C^k(N_\theta^+)}\right)
$$
as $t\to\infty$.
\item[(4)] There exists $Q_{\mathcal K_1,\mathcal K_2}'\in(L^1\cap L^\infty)(\ker\psi)$, depending only on $g,\mathcal K_1,\mathcal K_2$, such that
\begin{equation}
 |\err'(a_{tu+\sqrt t\,v},f,\phi)|\leq Q_{\mathcal K_1,\mathcal K_2}'(v)t^{-d'}\|f\|_{C^k(\mathcal M)}\|\phi\|_{C^k(N_\theta^+)}
\end{equation}
for all $t\geq1$ and $v\in\ker\psi$ with $a_{tu+\sqrt t\,v}\in A_\theta^+$.
\end{enumerate}
For fixed compact supports, the estimates in (3) and (4) hold uniformly as $g$ ranges over a fixed compact subset of $G$ satisfying $g\mathcal K_2.o\subset\widetilde{\mathcal M}$. The same uniformity holds in (2) for fixed $f,\phi$.
\end{proposition}
The proof of Proposition \ref{p1} follows the strategy of \cite{OhWinter}. 
We need the following lemma:
 
\begin{lemma}\label{fullmeas}
For any $g\in G$,
\begin{equation*}
 \nu (gN_\te^+(e^+) )>0\quad\text{and} \quad \nu_{ \op{i}}(gN_\te^-(e^-))>0.
\end{equation*}
\end{lemma}
\begin{proof}
Since $\Gamma$ is Zariski dense, $\La_{\Pi}$ is a Zariski dense subset of $\cal F_{\Pi}$ and the natural projections map $\Lambda_\Pi$ onto
$\Lambda_\theta$ and $\Lambda_{\i\theta}$ \cite{Ben}.
In particular, $\Lambda_\te$ and $\La_{\i\te}$ are Zariski dense in $\F_\te$ and $\F_{\i\te}$, respectively. 
Since $gN_\te^+(e^+)$ is a nonempty Zariski open subset of $\cal F_\te$, $gN_\te^+(e^+)\cap \Lambda_\te\neq\emptyset$. 
Similarly, $gN_\te^-(e^-)\cap\Lambda_{\i\te}\neq\emptyset$.
Since $\operatorname{supp}\nu=\Lambda_\te$ and $\operatorname{supp}\nu_\i=\Lambda_{\i\te}$, the lemma follows.
\end{proof}
We also need the following smoothness property 
\cite[Proposition 2.15]{OhShah}:
\begin{lemma}\label{PScont}
For any fixed $\rho\in C_c^\infty(N_\te^{\pm})$ and $g\in G$, the map $N_\te^{\mp}\to \br$ given by
$h \mapsto \mu_{gh, N_\te^{\pm}}^{\mathrm{PS}}(\rho)$ is $C^\infty$.
\end{lemma}
 
\begin{proof}
We prove the case $\rho\in C_c^\infty(N_\te^+)$; the other case is similar. For $h\in N_\te^-$ and
$\xi\in\mathcal F_\te$, define
$$
 \kappa(h,\xi)=
 \begin{cases}
 \rho(n)e^{\psi(\beta_\xi^\te(e,ghn))}, &\xi=(ghn)^+,\ n\in N_\te^+,\\
 0,&\text{otherwise}.
 \end{cases}
$$
The map $n\mapsto(ghn)^+$ is a diffeomorphism from $N_\te^+$ onto $ghN_\te^+(e^+)$. 
Its inverse depends smoothly on $(h,\xi)$. Compact support of $\rho$ makes its extension by zero smooth; all $h$-derivatives are uniformly bounded in $\xi$ when $h$ ranges over a compact set. 
Thus $\mu_{gh,N_\te^+}^{\PS}(\rho) =\int_{\mathcal F_\te}\kappa (h,\xi)\,d\nu(\xi)$ is smooth by differentiation under the integral. The same argument gives smooth dependence on $g\in G$.
\end{proof}
\begin{cor}\label{thickcor}
Given $\epsilon>0$ and $g\in G$, there exist $R>1$ and a nonnegative function $\rho_{g,\epsilon}\in
C_c^\infty((N_\te^-)_R)$ such that for all $ n\in (N_\te^+)_{\epsilon}$,
\begin{equation*}
 \mu_{gn,N_\te^-}^{\mathrm{PS}}(\rho_{g,\epsilon})>0.
\end{equation*}
\end{cor}
\begin{proof}
For each $j\in \NN$, let $\phi_j\in C_c^\infty((N_\te^-)_{j+1})$ be a nonnegative function such that $\phi_j|_{(N_\te^-)_j}=1$. 
By Lemma \ref{fullmeas}, for each $n\in N_\te^+$, there exists some $j_{n}\in \NN$ such that $\mu_{gn,N_\te^-}^{\mathrm{PS}}((N_\te^-)_{j_{n}})>0$. 
By Lemma \ref{PScont}, for each $n\in N_\te^+$, there exists $r_{n}>0$ such that
\begin{equation*}
 \mu_{gn_0,N_\te^-}^{\mathrm{PS}}(\phi_{j_n})>0 \quad \text{ for all } n_0\in \scrB(n):=\lbrace n_0\in N_\te^+\,:\, d(n,n_0)<r_{n}\rbrace.
\end{equation*}
By relative compactness, $(N_\te^+)_\epsilon\subset\cup_{i=1}^k \scrB(n_i)$ for some $n_1,\cdots,n_k\in N_\te^+$. 
Set
$R:=\max_i j_{n_i}+2$, and $\rho_{g,\epsilon}:=\sum_{i=1}^k\phi_{j_{n_i}}$.
Then $\rho_{g,\epsilon}\in C_c^\infty((N_\te^-)_R)$, and at each $n\in(N_\te^+)_\epsilon$, $\mu_{gn,N_\te^-}^{\PS}(\phi_{j_{n_i}})>0$ for at least one $i$.
The corollary follows.
\end{proof}
\noindent
 
\noindent{\bf Proof of Proposition \ref{p1}.}
 We may assume $f,\phi\geq0$. We first prove (1).
Fix $v\in \op{ker}\psi$ and write $a_t=a_{tu+\sqrt t\,v}$, which belongs to $A_\te^+$ for all sufficiently large $t$.
Choose $\epsilon_1>0$ such that $\phi\in C_c^\infty((N_\te^+)_{\epsilon_1})$. 
By Corollary \ref{thickcor}, there exist $R>1$ and a nonnegative $\rho_{g,\epsilon_1}\in C_c^\infty((N_\te^-)_R)$ such that
\begin{equation*}
 \mu_{gn,N_\te^-}^{\mathrm{PS}}(\rho_{g,\epsilon_1})>0\quad \text{ for all }n\in (N_\te^+)_{\epsilon_1}.
\end{equation*}
Fix $0<\epsilon<1$.
By Lemma \ref{lem.bump}, there exists a nonnegative $q_\epsilon\in C_c^\infty((A_\theta)_\epsilon)$ with integral one and $\|q_\epsilon\|_{C^k(A_\te)}\ll_k \epsilon^{-k-\dim A_\theta}$.
Then
\begin{align}\label{fthickening}
 &\int_{N_\te^+} f([gn] a_t)\phi(n)\,d\mu_{g,N_\te^+}^{\mathrm{PS}}(n) \notag\\
 &=\int_{N_\te^+} \left(\int_{N_\te^-A_\te}f([gn]a_t ) \frac{\phi(n)\rho_{g,\epsilon_1}(h)q_{\epsilon}(a)}{\mu_{gn,N_\te^-}^{\mathrm{PS}}(\rho_{g,\epsilon_1})} \,da\,d\mu_{gn,N_\te^-}^{\mathrm{PS}}(h)\!\right)d\mu_{g,N_\te^+}^{\mathrm{PS}}(n).
\end{align}
Define $\Phi_\epsilon\in C_c^\infty(g(N_\te^+)_{\epsilon_1}(N_\te^-)_R(A_\te)_\epsilon.o)\subset C_c^\infty(g\cal O)$ by
$$
\Phi_\epsilon(g_0):=
\begin{cases}
 \frac{\phi(n)\rho_{g,\epsilon_1}(h)q_\epsilon(a)} {\mu_{gn,N_\te^-}^{\PS}(\rho_{g,\epsilon_1})}& \text{if $g_0=gnha.o$ and $n\in(N_\te^+)_{\epsilon_1}$,}\\
 0&\text{otherwise.}
\end{cases}
$$
 Lemma \ref{PScont} shows $\Phi_{\epsilon}$ is smooth.
The coordinates $(n,h,a)\mapsto gnha.o$ and the bound for $q_\epsilon$ give
\begin{equation}\label{eq.phie}
 \|\Phi_\epsilon\|_{C^k(g\cal O)} \ll_{\rho,g,\epsilon_1} \|\phi\|_{C^k(N_\te^+)}\epsilon^{-(k+\dim A_\te)}.
\end{equation} 
Let $F_\epsilon^\pm$ be the functions associated with $f$ in \eqref{eq.ffe}, with the smoothing parameter replaced by $C_0\epsilon$, where $C_0\geq1$ is a fixed constant chosen below.
We claim that, for  all sufficiently large $t$, the following holds for all $n\in\supp\phi$, $h\in\supp\rho_{g,\epsilon_1}$ and $a\in\supp q_\epsilon$ such that $\bigl((gnh)^+,(gnh)^-\bigr)\in\Lambda_\theta^{(2)}$,
\begin{equation}\label{fplus}
 F_\epsilon^-([gnha]a_t)\leq f([gn]a_t)\leq F_\epsilon^+([gnha]a_t).
\end{equation}

Fix a compact subset $\mathcal K\subset\mathcal M$ containing the supports of $f$ and $F_\epsilon^\pm$ for all sufficiently small $\epsilon$. By Lemma \ref{Lem.FE1}, it suffices to show that, for all sufficiently large $t$, we have 
\begin{equation}\label{dm} d_{\mathcal M}([gn]a_t,[gnha]a_t)\leq C\epsilon
\end{equation}
whenever $a_t\in A_\theta^+$ and at least one of $[gn]a_t$ and $[gnha]a_t$ belongs to $\mathcal K$.
Choose a compact subset $\mathcal K_0\subset G$ mapping onto $\mathcal K$. Suppose first that $[gn]a_t\in\mathcal K$. Choose $q\in\mathcal K_0$ with $[q]=[gn]a_t$ and write $\gamma gna_t=qm$, where $\gamma\in\Gamma$ and $m\in M_\theta$.
 Since $a$ commutes with $a_t$ and $m$,
$$
 \gamma gnhaa_t.o=q\bigl(ma_t^{-1}ha_tm^{-1}\bigr)a.o.
$$
Since $\gamma^{-1}q.o=gna_t.o$, Lemma \ref{lem.return} gives $\mu_\theta(\gamma^{-1})=\log a_t+O(1)$.
Since $\psi(\log a_t)=t\psi(u)$ and $\psi(u)>0$,   $\gamma^{-1}\to \infty$ as $t\to\infty$. 
Regularity and Lemma \ref{lem.return} give $\gap(\gamma^{-1})\to\infty$ and  $\left\|\operatorname{Ad}(ma_t^{-1})\big|_{\mathfrak n_\te^-}\right\|\to0.$
Since $h$ ranges over a fixed compact subset of $N_\te^-$, $ ma_t^{-1}ha_tm^{-1}\in(N_\te^-)_{o(1)}.$
If $\Gamma$ is linearly $\theta$-regular, we have 
$$
 \gap(\gamma^{-1})\geq c\|\log a_t\|-C\geq c't-C,\qquad \left\|\operatorname{Ad}(ma_t^{-1})\big|_{\mathfrak n_\te^-}\right\|\leq Ce^{-ct},
$$
and hence $ma_t^{-1}ha_tm^{-1}\in(N_\te^-)_{Ce^{-ct}}$.
If $[gnha]a_t\in\mathcal K$, apply the same argument using $gn=(gnha)(a^{-1}h^{-1}a)a^{-1}$.
Here $a^{-1}h^{-1}a$ ranges over a fixed compact subset of $N_\theta^-$, and $a^{-1}\in(A_\theta)_\epsilon$. This proves \eqref{dm}. 
Choosing $C_0$ at least as large as the constant in \eqref{dm}, Lemma \ref{Lem.FE1} gives \eqref{fplus}.
 
Choose $\zeta\in C^\infty(\mathcal F_\theta^{(2)})$ whose support is contained in $\{(g^+,g^-):g.o\in\widetilde{\mathcal M}\}$ and which is identically $1$ on $\Lambda_\theta^{(2)}$.
 We use the same notation for its $A_\theta$-invariant extension to $G/M_\theta$.
By \eqref{fthickening}, \eqref{fplus} and Lemma \ref{NPdecomp}, we obtain
\begin{align*}
 &\Psi(t)\int_{N_\te^+} f([gn] a_t )\phi(n)\,d\mu_{g,N_\te^+}^{\mathrm{PS}}(n) \\
 &\leq \Psi(t)\int_{N_\te^+} \left(\int_{N_\te^-A_\te}F_{\epsilon}^+([gnha] a_t ) \Phi_\epsilon(g nha.o)\,da\,d\mu_{gn,N_\te^-}^{\mathrm{PS}}(h)\right)\,d\mu_{g,N_\te^+}^{\mathrm{PS}}(n) \\
 &=\Psi(t)\int_{\widetilde{\mathcal M}}F_{\epsilon}^+([g_0] a_t ) \Phi_\epsilon(g_0)\,d\tilde{m}^{\mathrm{BMS}}(g_0) \\
 &=\Psi(t)\int_{\mathcal M} F_{\epsilon}^+([g_0] a_t )\Phi_\epsilon^*([g_0])\,dm^{\mathrm{BMS}}([g_0]),
\end{align*}
where $\Phi_\epsilon^*([g_0]):=\sum_{\gamma\in \Gamma} \Phi_\epsilon(\gamma g_0)\zeta(\gamma g_0)$.  
Note that $\Phi_\epsilon^*\in C_c^\infty(\mathcal M)$ and that
\begin{equation}\label{eq.phie2}
 \norm{\Phi_\epsilon^*}_{C^k(\cal M)}\ll_{\zeta} \norm{\Phi_\epsilon}_{C^k(g\cal O)}\ll_{\rho,g,\epsilon_1} \epsilon^{-(k+\op{dim}A_\te)}\norm{\phi}_{C^k(N_\te^+)}
\end{equation}
by \eqref{eq.phie}.
Using the definition of $\err$,
\begin{align*}
 &\Psi(t)\int_{N_\te^+} f([gn] a_t )\phi(n)\,d\mu_{g,N_\te^+}^{\mathrm{PS}}(n)\leq J(v) m^{\mathrm{BMS}}(F_{\epsilon}^+)m^{\mathrm{BMS}}(\Phi_\epsilon^*)+\mathsf{E}(a_t,F_{\epsilon}^+,\Phi_\epsilon^*)\\
 &=J(v) m^{\mathrm{BMS}}(F_{\epsilon}^+)\tilde{m}^{\mathrm{BMS}}( \Phi_\epsilon\zeta)+\mathsf{E}(a_t,F_{\epsilon}^+,\Phi_\epsilon^*).
\end{align*}
Note that $\tilde{m}^{\mathrm{BMS}}( \Phi_\epsilon\zeta)=\tilde{m}^{\mathrm{BMS}}( \Phi_\epsilon)$. 
Using Lemma \ref{NPdecomp}, we have
\begin{align*}
 &\tilde{m}^{\mathrm{BMS}}( \Phi_\epsilon)=\int_{N_\te^+} \left(\int_{N_\te^-A_\te}\Phi_\epsilon(g nha.o)\,da\,d\mu_{gn,N_\te^-}^{\mathrm{PS}}(h)\right)\,d\mu_{g,N_\te^+}^{\mathrm{PS}}(n) \\
 &=\int_{N_\te^+} \frac{\phi(n)}{\mu_{gn,N_\te^-}^{\mathrm{PS}}(\rho_{g,\epsilon_1})}\left(\int_{N_\te^-A_\te} \rho_{g,\epsilon_1}(h)q_{\epsilon}(a)\,da\,d\mu_{gn,N_\te^-}^{\mathrm{PS}}(h)\!\right)\,d\mu_{g,N_\te^+}^{\mathrm{PS}}(n) \\
 &= \mu_{g,N_\te^+}^{\mathrm{PS}}(\phi).
\end{align*}
Hence
\begin{equation}\label{eq.alt}
 \Psi(t)\int_{N_\te^+} f([gn] a_t )\phi(n)\,d\mu_{g,N_\te^+}^{\mathrm{PS}}(n)\leq J(v)\, m^{\mathrm{BMS}}(f)\,\mu_{g,N_\te^+}^{\mathrm{PS}}(\phi)+\err_{\epsilon}^+(a_t,f,\phi)
\end{equation}
where
$\err_{\epsilon}^\pm(a_t,f,\phi):=J(v)\, m^{\mathrm{BMS}}(F_{\epsilon}^\pm-f)\,\mu_{g,N_\te^+}^{\mathrm{PS}}(\phi)+\mathsf{E}(a_t,F_{\epsilon}^\pm,\Phi_\epsilon^*)$.
 
Similarly, replacing $F^+_{\epsilon}$ with $F^-_{\epsilon}$ gives the lower bound
\begin{equation}\label{eq.alt2}
 \Psi(t)\int_{N_\te^+} f([gn] a_t )\phi(n)\,d\mu_{g,N_\te^+}^{\mathrm{PS}}(n)\geq J(v)\, m^{\mathrm{BMS}}(f)\,\mu_{g,N_\te^+}^{\mathrm{PS}}(\phi)+\err_{\epsilon}^-(a_t,f,\phi).
\end{equation}
 Letting $t\to\infty$ first and then
$\epsilon\to0$ in \eqref{eq.alt} and \eqref{eq.alt2}, Lemma \ref{lem.e1} gives
Proposition \ref{p1}(1).
 
 For (2), assume condition (2) of Definition \ref{lmp}. By Lemma \ref{lem.return}, $ma_t^{-1}ha_tm^{-1}$ lies in a fixed compact subset of $N_\te^-$ whenever $a_t\in A_\te^+$.
Right multiplication by $N_\te^-$ preserves the forward endpoint, which belongs to $\Lambda_\te$ here. 
So when $[gn]a_t\in\supp f$, all the points $[gnha]a_t$ in \eqref{fthickening} lie in a fixed compact subset of $\mathcal M$. Choose $B\in C_c^\infty(\mathcal M)$, $B\geq0$, with $B\geq\|f\|_\infty$ on that set. 
The preceding calculation gives
$$
 \Psi(t)\int_{N_\te^+}f([gn]a_t)\phi(n)\,d\mu_{g,N_\te^+}^{\PS}(n) \leq\Psi(t)\int_{\mathcal M}B(xa_t)\Phi_\epsilon^*(x) \,dm^{\BMS}(x).
$$
Definition \ref{lmp}(2) bounds the right-hand side.  Choosing fixed smooth functions equal to one on these compact supports gives $C'\ll\|f\|_\infty\|\phi\|_\infty$.
 
We next prove (3) and (4). By Lemmas \ref{Lem.FE1} and \ref{lem.e1},
$$
 \bigl|m^{\BMS}(F_{\epsilon}^\pm-f)\bigr| \ll\epsilon\|f\|_{C^1},\qquad \|F_{\epsilon}^\pm\|_{C^k} \ll\epsilon^{-k}\|f\|_{C^k}.
$$
Let $\mathcal K\subset\mathcal M$ be a fixed compact set containing the supports of $F_{\epsilon}^\pm$ and $\Phi_\epsilon^*$ for all sufficiently small $\epsilon$.
Local finiteness also gives $\mu_{g,N_\te^+}^{\PS}(\phi)\ll\|\phi\|_\infty$. For (3), fix $v$ and use \eqref{eq.kd}, \eqref{eq.phie2}, \eqref{eq.alt} and \eqref{eq.alt2} to obtain
$$
 |\err'(a_t,f,\phi)|\ll_{g,\mathcal K_1,\mathcal K_2,v}\left(\epsilon+t^{-d}\epsilon^{-2k-\dim A_\te}\right)\|f\|_{C^k}\|\phi\|_{C^k},
$$
provided $e^{-ct}\ll\epsilon$ and $t$ is sufficiently large depending on the compact supports and $v$. Taking $\epsilon$ to be a fixed sufficiently small multiple of $t^{-d/(2k+\dim A_\te+1)}$ proves (3), with $d'=d/(2k+\dim A_\te+1)$.
 
For (4), the same calculation using \eqref{eq.kd2} gives
$$
 |\err'(a_t,f,\phi)| \ll\|f\|_{C^k}\|\phi\|_{C^k} \left(\epsilon J(v) +\epsilon^{-2k-\dim A_\te}t^{-d}Q_{\mathcal K}(v)\right),
$$
provided $e^{-ct}\ll\epsilon$. With the same choice of $\epsilon$, this proves (4) for all  $t\ge T_0$, with the bound $C(J+Q_{\mathcal K})$, where $T_0$ is independent of $v$.
For $1\leq t\leq T_0$, argue as in the proof of (2) with $\epsilon$ fixed. Applying \eqref{eq.kd2} to the smooth functions $B$ and $\Phi_\epsilon^*$ used there, and enlarging $\mathcal K$ to contain their supports, gives
$$
 \Psi(t)\int_{N_\te^+} f([gn]a_t)\phi(n)\,d\mu_{g,N_\te^+}^{\PS}(n)\ll\bigl(J(v)+Q_{\mathcal K}(v)\bigr)\|f\|_{C^k}\|\phi\|_{C^k}.
$$
The main term has the same bound.
Since $t^{-d'}\geq T_0^{-d'}$ for $t\leq T_0$, this proves (4) for every $t\geq1$, with $Q'_{\mathcal K_1,\mathcal K_2}(v)=C\bigl(J(v)+Q_{\mathcal K}(v)\bigr)\in(L^1\cap L^\infty)(\ker\psi)$.
The constants can be chosen uniformly for $g$ in a fixed compact subset of $G$ and $\supp\phi\subset\mathcal K_2$, provided $g\mathcal K_2.o\subset\widetilde{\mathcal M}$ for every $g$ in that set. \qed
 
\section{Translates of Lebesgue measures and Haar local mixing}\label{s2}
We keep the notation and assumptions of Section \ref{sec.set}.
We first prove the analogue of Proposition \ref{p1} with Lebesgue measure $dn$ in place of $d\mu_{g,N_\te^+}^{\PS}$ (see \eqref{eq.Leb2}), and then deduce local mixing for $dx$.
 
\begin{proposition}\label{p2} \rm
Under conditions (1) and (2) of Definition \ref{lmp}, (1) and (2) below hold, respectively.
\begin{enumerate}
\item For all $g\in G$, $v\in\ker\psi$, $f\in C_c^\infty(\mathcal M)$, and $\phi\in C_c^\infty(N_\te^+)$ with $g(\supp\phi).o\subset\widetilde{\mathcal M}$,
\begin{equation*}
\begin{split}
 &\Psi(t) e^{(2\rho_\te-\psi)(\log a_{tu+\sqrt t\,v})} \int_{N_\te^+}f\big([gn] \,a_{tu+\sqrt t\,v}\big)\phi(n)\,dn\\
 &=J(v)\,m^{\mathrm{BR}}(f)\,\mu_{g,N_\te^+}^{\mathrm{PS}}(\phi)+\err''(a_{tu+\sqrt t\,v},f,\phi),
\end{split}
\end{equation*}
and $\err''(a_{tu+\sqrt t\,v},f,\phi)\to 0$ as $t\to\infty$.
\item There exists $C''=C''(g,f,\phi)>0$ such that for all $t\geq1$ and $v\in\ker\psi$ with $a_{tu+\sqrt t\,v}\in A_\te^+$,
\begin{equation*}
 \left| \Psi(t)e^{ (2\rho_\te-\psi)(\log a_{tu+\sqrt t\,v})}\int_{N_\te^+}f\big([gn] \,a_{tu+\sqrt t\,v}\big)\phi(n)\,dn\right|\leq C''.
\end{equation*}
\end{enumerate}
Assume now that $\Gamma$ is linearly $\theta$-regular \eqref{eq.rc}. Then (3) and (4) below hold under conditions (1) and (2) of Definition \ref{def.em}, respectively, for some $d''>0$ depending only on $d$, $k$ and $G$.
\begin{enumerate}
\item[(3)] For every $v\in\ker\psi$, compact sets $\mathcal K_1\subset\mathcal M$, $\mathcal K_2\subset N_\te^+$ with $g\mathcal K_2.o\subset\widetilde{\mathcal M}$, and $f\in C_c^\infty(\mathcal K_1)$, $\phi\in C_c^\infty(\mathcal K_2)$, we have
$$
 \err''(a_{tu+\sqrt t\,v},f,\phi)= O_{g,\cal K_1,\cal K_2,v}(t^{-d''}\norm{f}_{C^k(\mathcal M)}\norm{\phi}_{C^{k}(N_\te^+)}) \qquad \text{as $t\to\infty$}.
$$
\item[(4)] For all $(t,v)$ with $t\geq1$ and $a_{tu+\sqrt t\,v}\in A_\te^+$,
\begin{equation*}
 |\err''(a_{tu+\sqrt t\,v},f,\phi)|\leq Q_{\cal K_1,\cal K_2}''(v) t^{-d''}\norm{f}_{C^k(\mathcal M)}\norm{\phi}_{C^k(N_\te^+)}
\end{equation*}
for some $Q_{\cal K_1,\cal K_2}''\in (L^1\cap L^\infty)(\ker\psi)$, depending only on $g,\mathcal K_1,\mathcal K_2$.
\end{enumerate}
The constants in (3) and the lower bound on $t$ are uniform for $f\in C_c^\infty(\mathcal K_1)$ and $\phi\in C_c^\infty(\mathcal K_2)$. 
They can also be chosen uniformly for $g$ in a fixed compact subset of $G$ satisfying $g\mathcal K_2.o\subset\widetilde{\mathcal M}$; the same holds for the constant in (2) and the function in (4).
\end{proposition}
\noindent
For $0<\epsilon<1$, let $\tilde \eta_{\epsilon}\in C_c^\infty((N_\te^+)_{\epsilon})$ be as in Lemma \ref{lem.bump}, satisfying
\begin{equation}\label{eq.eta2}
 \int_{N_\te^+} \tilde \eta_{\epsilon}(h)\,dh=1,\quad\text{ and }\quad \norm{\tilde \eta_{\epsilon}}_{C^k(N_\te^+)}\ll_k {\epsilon}^{-(k+\op{dim}N_\te^+)}.
\end{equation}
For $\phi\in C(N_\te^+)$ and $\epsilon>0$, define $\phi_{\epsilon}^\pm\in C(N_\te^+)$ by
\begin{equation*}
 \phi^+_{\epsilon}(n):=\sup_{b\in (N_\te^+)_{\epsilon}} \phi(nb),\quad\text{ and }\quad \phi^-_{\epsilon}(n):=\inf_{b\in (N_\te^+)_{\epsilon}} \phi(nb).
\end{equation*}
For $\phi\in C(N_\te^+)$, its convolution with $\tilde \eta_\epsilon$ is defined by
$$
 (\phi*\tilde\eta_\epsilon)(n) =\int_{N_\te^+}\phi(nb)\tilde\eta_\epsilon(b^{-1})\,db.
$$
\begin{lemma}\label{lem.pp2}
For $\phi\in C(N_\te^+)$, $n\in N_\te^+$, $n'\in (N_\te^+)_{\epsilon}$ and
$0<\epsilon<1$, we have
$$
 (\phi_{2\epsilon}^-*\tilde \eta_{\epsilon})(n)\leq \phi(nn')\leq (\phi_{2\epsilon}^+*\tilde \eta_{\epsilon})(n).
$$
\end{lemma}
\begin{proof}
For $b,n'\in(N_\te^+)_\epsilon$, we have $(nb)^{-1}(nn')=b^{-1}n'\in(N_\te^+)_{2\epsilon}$. 
Hence
$$
 \phi_{2\epsilon}^+*\tilde \eta_{\epsilon}(n)= \int_{N_\te^+}\phi_{2\epsilon}^+(nb)\tilde\eta_{\epsilon}(b^{-1})\,db\geq \int_{N_\te^+}\phi(nn')\tilde\eta_{\epsilon}(b^{-1})\,db=\phi(nn').
$$
The lower bound is obtained similarly.
\end{proof}
\begin{lemma}\label{lem.pp3}
Let $\phi\in C^1(N_\te^+)$ and $0<\epsilon<1$. We have:
\begin{enumerate}
\item for all $n\in N_\te^+$,
$\phi_{\epsilon}^+(n)\leq \phi_{3\epsilon}^+*\tilde\eta_\epsilon(n) \text{ and
}\phi_{\epsilon}^-(n)\geq \phi_{3\epsilon}^-*\tilde\eta_\epsilon(n)$.
\item for all $k$, $\norm{\phi_{3\epsilon}^\pm*\tilde\eta_{\epsilon}}_{C^k(N_\te^+)}\ll_k \epsilon^{-(k+\op{dim}N_\te^+)}\norm{\phi}_{\infty}$.
\item $\norm{\phi_{3\epsilon}^+-\phi_{3\epsilon}^-}_{\infty}\ll \epsilon\norm{\phi}_{C^1(N_\te^+)}$.
\end{enumerate}
\end{lemma}
\begin{proof}
(1) follows from Lemma \ref{lem.pp2} and the inclusion $(N_\te^+)_\epsilon(N_\te^+)_{2\epsilon}\subset(N_\te^+)_{3\epsilon}$.
For (2), differentiating the convolution, with the derivatives on $\tilde\eta_\epsilon$, gives
$$
 \| \phi_{3\epsilon}^\pm*\tilde\eta_\epsilon\|_{C^k}\ll_k\|\phi_{3\epsilon}^\pm\|_\infty\|\tilde\eta_\epsilon\|_{C^k}\ll_k\epsilon^{-(k+\dim N_\te^+)}\|\phi\|_{\infty}.
$$
Here the first bound uses that $\supp\tilde\eta_\epsilon$ lies in the fixed compact set
$\overline{(N_\te^+)_1}$, and the second uses \eqref{eq.eta2}.
(3) is the mean value estimate.
\end{proof}
 
\noindent\textbf{Proof of Proposition \ref{p2}.}
We first prove (1).
For $g\in G$ and $\epsilon>0$, set
$ \scrB_\epsilon(g) :=g(N_\te^-)_\epsilon(A_\te)_\epsilon(N_\te^+)_\epsilon.o.$
For $g.o\in\widetilde{\mathcal M}$, choose $\epsilon_0(g)>0$ such that $\scrB_{\epsilon_0(g)}(g)\subset\widetilde{\mathcal M}$ and its projection to $\mathcal M$ is
injective. Fix $v\in\ker\psi$ and write $a_t=a_{tu+\sqrt t\,v}$. 
By linearity and smooth partitions of unity, it suffices to consider $f,\phi\geq0$ with
$$
 f\in C_c^\infty(\Gamma\scrB_\epsilon(g_0)),\qquad \phi\in C_c^\infty((N_\te^+)_\epsilon),\qquad 0<\epsilon<\min\{1,\epsilon_0(g_0)\}.
$$
Choose $\tilde f\in C_c^\infty(\scrB_\epsilon(g_0))$, $\tilde f\geq0$, such that for $q.o\in\widetilde{\mathcal M}$,
\begin{equation*}
 f([q])=\sum_{\gamma\in\Gamma}\tilde f(\gamma q.o).
\end{equation*}
If $\tilde f(\ga gna_t.o)\phi(n)\neq 0$, then $\ga gna_t \in g_0(N_\te^-)_\epsilon(A_\te)_\epsilon(N_\te^+)_\epsilon M_\te$ and hence
$$
\gamma\in  \Ga_{t,v}:=\{\ga\in\Ga:\gamma ga_t\in g_0(N_\te^-)_\epsilon(A_\te)_\epsilon M_\te N_\te^+\}.
$$
For each $\gamma\in\Ga_{t,v}$, write 
$
 \gamma ga_t=g_0h_\gamma a_\gamma m_\gamma n_\gamma \in  
 g_0(N_\te^-)_\epsilon(A_\te)_\epsilon M_\te N_\te^+.
$
Then
\begin{align}\label{eq.un1}
&\int_{N_\te^+}f([gn]a_t)\phi(n)\,dn 
=\sum_{\gamma\in\Gamma_{t,v}} \int_{N_\te^+}\tilde f(\gamma ga_t(a_t^{-1}na_t).o)\phi(n)\,dn\notag \\
&=e^{-2\rho_\te(\log a_t)}\sum_{\gamma\in\Gamma_{t,v}} \int_{N_\te^+}\tilde f(\gamma ga_tn.o)\phi(a_tna_t^{-1})\,dn \notag \\
&=e^{-2\rho_\te(\log a_t)}\sum_{\gamma\in\Gamma_{t,v}} \int_{N_\te^+}\tilde f(g_0 h_\ga a_\ga m_\ga n_\ga n.o)\phi(a_tna_t^{-1})\,dn \notag \\
&=e^{-2\rho_\te(\log a_t)}\sum_{\gamma\in\Gamma_{t,v}} \int_{N_\te^+}\tilde f(g_0 h_\ga a_\ga   n.o)\phi(a_tn_\ga^{-1}m_\ga^{-1}n m_\ga a_t^{-1})\,dn 
\end{align}
where the last step follows from a change of variable and the fact that $|\det(\operatorname{Ad}(m)|_{\mathfrak n_\te^+})|=1$ for all $m\in M_\te$.
For $\ga\in\Ga_{t,v}$ and $z\in N_\te^+$, define 
\begin{equation*}
\mathsf{n}_{\gamma,t}(z) :=a_tn_\gamma^{-1}m_\gamma^{-1}zm_\gamma a_t^{-1}\in N_\te^+.    
\end{equation*}
In particular, for $z\in N_\te^+$ we have
\begin{equation}\label{eq.x3}
\gamma g\mathsf n_{\gamma,t}(z)a_t
=g_0h_\gamma a_\gamma z m_\gamma.
\end{equation}
Then \eqref{eq.un1} becomes
\begin{equation}\label{intchain}
 e^{2\rho_\te(\log a_t)} \int_{N_\te^+}f([gn]a_t)\phi(n)\,dn =\sum_{\gamma\in\Gamma_{t,v}}\int_{N_\te^+} \tilde f(g_0h_\gamma a_\gamma z.o) \phi(\mathsf{n}_{\gamma,t}(z))\,dz.
\end{equation}
By the proof of Corollary \ref{thickcor} with $N_\te^+$ and $N_\te^-$ exchanged, there exist $R>0$ and a nonnegative $\rho\in C_c^\infty((N_\te^+)_R)$ such that $\mu_{g_0ha,N_\te^+}^{\PS}(\rho)>0$ for all $ha\in (N_\te^-)_{\epsilon}(A_\te)_{\epsilon}$. 
 
Define $\tilde F \in C_c^\infty(g_0 (N_\te^-)_\epsilon (A_\te)_{\epsilon} (N_\te^+)_R.o)$ by
\begin{equation*}
\tilde F(x) :=
\begin{cases}
\frac{\rho(z)}{\mu_{g_0ha,N_\te^+}^{\PS}(\rho)} \int_{N_\te^+}\tilde f(g_0haw.o)\,dw&\text{if $x=g_0haz.o\in g_0 (N_\te^-)_\epsilon (A_\te)_{\epsilon} (N_\te^+)_R.o$}\\
0&\text{otherwise.}
\end{cases}  
\end{equation*}
We have, for all $ha\in (N_\te^-)_{\epsilon}(A_\te)_{\epsilon}$,
\begin{equation}\label{eq.tilf2}
 \int_{N_\te^+}\tilde F(g_0haz.o) \,d\mu_{g_0ha,N_\te^+}^{\PS}(z)=\int_{(N_\te^+)_{\epsilon}}\tilde f(g_0haz.o)\,dz.
\end{equation}
Choose $\zeta\in C_c^\infty(\widetilde{\cal M})$ which is identically one on the compact set $\supp\tilde F\cap\{q.o:q^+\in\Lambda_\te\}\subset\widetilde{\cal M}$.
Define $F\in C_c^\infty(\mathcal M)$ by:
\begin{equation*}
 F([q]):=\sum_{\gamma\in\Gamma}(\zeta\tilde F)(\gamma q.o).
\end{equation*}
 
Choose a compact set $S\subset G$ with $S.o\subset\widetilde{\mathcal M}$ containing $g\supp\phi$ and the elements $g_0haz$ with $g_0haz.o$ in $\supp\tilde f$ or in $\supp\tilde F\cap\{q.o:q^+\in\Lambda_\theta\}$.
Fix a compact set $D\subset N_\te^+$ containing $\overline{(N_\te^+)_\epsilon}\cup\supp\rho$.
The following claim is a consequence of Lemma \ref{lem.return}:
 
\textbf{Claim.} There exist $T_0>0$ and a function $\omega: [T_0,\infty)\to (0,\infty)$ with $\omega(t)\to 0$ as $t\to\infty$ such that the following holds: 
whenever $t\geq T_0$, $v\in\ker\psi$, $a_t\in A_\te^+$ and
$\gamma\in\Gamma_{t,v}$ satisfy
$\gamma p a_t=qm_\gamma$ for some $p,q\in S$, we have 
$$
 \mathsf{n}_{\gamma,t}(z)^{-1}\mathsf{n}_{\gamma,t}(w)\in(N_\te^+)_{\omega(t)}
$$
for all $z$, $w\in D$.
If $\Gamma$ is linearly $\theta$-regular, then $\omega$ can be chosen as $\omega(t)=C_0e^{-ct}$ for some $C_0,c>0$.
The choices are independent of $v,\gamma,z,w$.
 
Indeed,
$$
 \mathsf n_{\gamma,t}(z)^{-1}\mathsf n_{\gamma,t}(w)=a_tm_\gamma^{-1}(z^{-1}w)m_\gamma a_t^{-1}.
$$
The second part of Lemma \ref{lem.return}, applied to $\gamma p a_t=qm_\gamma$, gives $\mu_\theta(\gamma^{-1})=\log a_t+O_S(1)$. Since $\|\log a_t\|\geq t\psi(u)/\|\psi\|$, we have $\gamma^{-1}\to\infty$ as $t\to\infty$, uniformly in $v$. Regularity and \eqref{eq.holo} therefore prove the claim, since $D^{-1}D$ is compact. Under linear $\theta$-regularity, $\gap(\gamma^{-1})\geq ct-C$, giving the stated exponential bound.
 
We now prove that, for every $t\geq T_0$, $v\in\ker\psi$ with $a_t\in A_\te^+$, every $\epsilon_1>0$ satisfying $\omega(t)\leq\epsilon_1$, and every $\gamma\in\Gamma_{t,v}$, the following inequalities hold:
\begin{align}
 &\int_{N_\te^+}\tilde F(g_0h_\gamma a_\gamma w.o) \phi^-_{\epsilon_1}(\mathsf{n}_{\gamma,t}(w)) \,d\mu_{g_0h_\gamma a_\gamma,N_\te^+}^{\PS}(w)\label{eq.u1}\\
 &\quad\leq\int_{N_\te^+}\tilde f(g_0h_\gamma a_\gamma z.o) \phi(\mathsf{n}_{\gamma,t}(z))\,dz\label{eq.u2}\\
 &\quad\leq\int_{N_\te^+}\tilde F(g_0h_\gamma a_\gamma w.o) \phi^+_{\epsilon_1}(\mathsf{n}_{\gamma,t}(w)) \,d\mu_{g_0h_\gamma a_\gamma,N_\te^+}^{\PS}(w).\notag
\end{align}
For the upper inequality, there is nothing to prove if the integral in \eqref{eq.u2} is zero. Otherwise, choose $z_0$ such that
$$
 \tilde f(g_0h_\gamma a_\gamma z_0.o)>0,\qquad
 \phi(\mathsf n_{\gamma,t}(z_0))>0.
$$
Note that $g\mathsf n_{\gamma,t}(z_0)\in S$, $g_0h_\gamma a_\gamma z_0 \in S$, and
$$
 \gamma (g\mathsf n_{\gamma,t}(z_0))a_t
 =(g_0h_\gamma a_\gamma z_0)m_\gamma
$$
by \eqref{eq.x3}.
The claim and $\omega(t)\leq\epsilon_1$ imply
$$
 \sup_{z\in(N_\te^+)_\epsilon}
 \phi(\mathsf n_{\gamma,t}(z))
 \leq
 \phi^+_{\epsilon_1}(\mathsf n_{\gamma,t}(w))
$$
for all $w\in\supp\rho$.
Hence
\begin{align*}
 &\int_{N_\te^+}
   \tilde f(g_0h_\gamma a_\gamma z.o)
   \phi(\mathsf n_{\gamma,t}(z))\,dz\leq
   \left(
    \sup_{z\in(N_\te^+)_\epsilon}
    \phi(\mathsf n_{\gamma,t}(z))
   \right)
   \int_{N_\te^+}\tilde f(g_0h_\gamma a_\gamma z.o)\,dz\\
 &=
   \int_{N_\te^+}\tilde F(g_0h_\gamma a_\gamma w.o)
   \left(
    \sup_{z\in(N_\te^+)_\epsilon}
    \phi(\mathsf n_{\gamma,t}(z))
   \right)
   \,d\mu_{g_0h_\gamma a_\gamma,N_\te^+}^{\PS}(w)\\
 &\leq
   \int_{N_\te^+}\tilde F(g_0h_\gamma a_\gamma w.o)
   \phi^+_{\epsilon_1}(\mathsf n_{\gamma,t}(w))
   \,d\mu_{g_0h_\gamma a_\gamma,N_\te^+}^{\PS}(w).
\end{align*}
Here the equality follows from \eqref{eq.tilf2} and the condition on $\supp \tilde f$.
 
For the lower inequality, there is nothing to prove if the integral in \eqref{eq.u1} is zero. Otherwise, there exists $w_0$ such that
$$
\tilde F(g_0h_\gamma a_\gamma w_0.o)>0,\quad
\phi^-_{\epsilon_1}(\mathsf n_{\gamma,t}(w_0))>0,\quad
g_0h_\gamma a_\gamma w_0^+\in\Lambda_\te.
$$

Then $g\mathsf n_{\gamma,t}(w_0)$ and $g_0h_\gamma a_\gamma w_0$ belong to $S$, and by \eqref{eq.x3} the claim applies. The claim implies that for all $w\in\supp\rho$,
$$
 \phi^-_{\epsilon_1}(\mathsf n_{\gamma,t}(w))\leq\inf_{z\in(N_\theta^+)_\epsilon}\phi(\mathsf n_{\gamma,t}(z)).
$$
Integrating and using \eqref{eq.tilf2} proves the lower inequality.
 
Choose a compact set $Q\subset N_\te^+$ containing the supports of $\phi^\pm_{3\delta}*\tilde\eta_\delta$ for $0<\delta<1$. The compact set $\{n\in
Q:gn^+\in\Lambda_\te\}$ is contained in the open set $\{n:gn.o\in\widetilde{\mathcal M}\}$. 
Choose $\chi_0\in C_c^\infty(N_\te^+)$ supported in this open set, with $0\leq\chi_0\leq1$ and equal to one on that compact set, and define
$$
 \tilde\phi^\pm_\delta :=\chi_0(\phi^\pm_{3\delta}*\tilde\eta_\delta).
$$
Multiplication by $\chi_0$ leaves all integrals against $\mu_{g,N_\te^+}^{\PS}$ unchanged. 
By Lemma \ref{MANconj} and the identity $\mu_{\gamma g,N_\te^+}^{\PS}=\mu_{g,N_\te^+}^{\PS}$ for $\gamma\in\Gamma$,
\begin{equation}\label{eq.x2}
 d\mu_{g,N_\te^+}^{\PS}(\mathsf{n}_{\gamma,t}(z)) =e^{-\psi(\log a_t)} d\mu_{g_0h_\gamma a_\gamma,N_\te^+}^{\PS}(z).
\end{equation}
Summing the preceding inequalities and using \eqref{eq.x2}, \eqref{intchain}, and Lemma \ref{lem.pp3}(1), we obtain
\begin{align}\label{eq.psi}
 &\Psi(t)\int_{N_\te^+}F([gn]a_t) \tilde\phi^-_{\epsilon_1}(n)\,d\mu_{g,N_\te^+}^{\PS}(n)\notag\\
 &\quad\leq\Psi(t)e^{(2\rho_\te-\psi)(\log a_t)} \int_{N_\te^+}f([gn]a_t)\phi(n)\,dn\\
 &\quad\leq\Psi(t)\int_{N_\te^+}F([gn]a_t) \tilde\phi^+_{\epsilon_1}(n)\,d\mu_{g,N_\te^+}^{\PS}(n).\notag
\end{align}
For fixed $v$, we have $a_t\in A_\te^+$ for all sufficiently large $t$. Letting $t\to\infty$ and then $\epsilon_1\to0$ in \eqref{eq.psi}, Proposition \ref{p1}(1) gives the limit $J(v)m^{\BMS}(F)\mu_{g,N_\te^+}^{\PS}(\phi)$. Moreover, \eqref{PNgood}, \eqref{eq.tilf2}, and \eqref{BRPN} give
\begin{align*}
 m^{\BMS}(F) &=\int_{N_\te^-A_\te}\int_{N_\te^+} \tilde F(g_0haz.o)\,d\mu_{g_0ha,N_\te^+}^{\PS}(z) \,d\mu_{g_0,N_\te^-A_\te}^{\PS}(ha)\\
 &=\int_{N_\te^-A_\te}\int_{N_\te^+} \tilde f(g_0haz.o)\,dz \,d\mu_{g_0,N_\te^-A_\te}^{\PS}(ha)=m^{\BR}(f).
\end{align*}
This proves (1).
 
For (2), assume condition (2) of Definition \ref{lmp}. By Lemma \ref{lem.return}, the upper inequality following \eqref{eq.u2} holds with $\epsilon_1$ replaced by a fixed $\delta_0>0$ whenever $a_t\in A_\te^+$.
Choose $B\in C_c^\infty(N_\te^+)$ with $g(\supp B).o\subset\widetilde{\mathcal M}$ and $B\geq\phi_{\delta_0}^+$ almost everywhere for $\mu_{g,N_\te^+}^{\PS}$. Such a function exists by the argument used for $\chi_0$.
Hence, when $a_t\in A_\te^+$,
$$
\begin{aligned}
 &\Psi(t)e^{(2\rho_\te-\psi)(\log a_t)} \int_{N_\te^+}f([gn]a_t)\phi(n)\,dn\\
 &\leq\Psi(t)\int_{N_\te^+}F([gn]a_t)B(n)\,d\mu_{g,N_\te^+}^{\PS}(n).
\end{aligned}
$$
Proposition \ref{p1}(2) bounds the right-hand side independently of $(t,v)$.
 
For (3) and (4), assume linear $\theta$-regularity and the corresponding condition of Definition \ref{def.em}. The definitions and Lemma \ref{lem.pp3} give
\begin{gather*}
 \|F\|_{C^k}\ll\|f\|_{C^k},\qquad \|\tilde\phi^\pm_{\epsilon_1}\|_{C^k}\ll\epsilon_1^{-k-\dim N_\te^+}\|\phi\|_{C^k},\\
 \left|\mu_{g,N_\te^+}^{\PS}(\tilde\phi^\pm_{\epsilon_1}-\phi)\right|\ll\epsilon_1\|\phi\|_{C^1}.
\end{gather*}
For (3), fix $v$. Proposition \ref{p1}(3) and \eqref{eq.psi} give
$$
 |\err''(a_t,f,\phi)|\ll_{g,\mathcal K_1,\mathcal K_2,v}\left(\epsilon_1+t^{-d'}\epsilon_1^{-k-\dim N_\te^+}\right)\|f\|_{C^k}\|\phi\|_{C^k}.
$$
Taking $\epsilon_1$ to be a fixed sufficiently small multiple of $t^{-b}$, where $b=d'/(k+\dim N_\te^++1)$, proves (3) with $d''=b$.
 
For (4), Proposition \ref{p1}(4) and \eqref{eq.psi} yield
$$
 |\err''(a_t,f,\phi)| \ll\|f\|_{C^k}\|\phi\|_{C^k} \left(\epsilon_1J(v) +t^{-d'}\epsilon_1^{-k-\dim N_\te^+}Q'(v)\right).
$$
The same choice of $\epsilon_1$ proves (4) for all sufficiently large $t$, with $d''=b$ and the bound $C(J+Q')$. For $1\leq t\leq T_0$, use the inequality in the proof of (2), followed by Proposition \ref{p1}(4). This bounds the normalized integral and its main term by
$$
 C\bigl(J(v)+Q'(v)\bigr)\|f\|_{C^k}\|\phi\|_{C^k}.
$$
Since $t^{-d''}\geq T_0^{-d''}$ for $t\leq T_0$, we may take
$
 Q''(v)=C\bigl(J(v)+Q'(v)\bigr)\in(L^1\cap L^\infty)(\ker\psi).
$
Here $J$ is bounded by assumption and integrable by Definition \ref{def.em}(2). This proves (4).
The constants, compact sets, and $Q''$ can be chosen uniformly for $g$ in a fixed compact subset of $G$ and $\supp\phi\subset\mathcal K_2$, provided $g\mathcal K_2.o\subset\widetilde{\mathcal M}$ for every $g$ in that set. \qed
 
We now deduce local mixing for the Haar measure.

\begin{theorem}\label{prop.mixH0}
Under conditions (1) and (2) of Definition \ref{lmp}, (1) and (2) below hold, respectively.
\begin{enumerate}
\item For every $v\in\ker\psi$ and $f_1,f_2\in C_c^\infty(\mathcal M)$,
$$
\begin{aligned}
 &\Psi(t)e^{(2\rho_\te-\psi)(\log a_{tu+\sqrt t\,v})}\int_{\mathcal M}f_1(xa_{tu+\sqrt t\,v})f_2(x)\,dx\\
 &=J(v)m^{\BR}(f_1)m^{\BR_*}(f_2)+\err'''(a_{tu+\sqrt t\,v},f_1,f_2),
\end{aligned}
$$
where $\err'''(a_{tu+\sqrt t\,v},f_1,f_2)\to0$ as $t\to\infty$.
\item For each $f_1,f_2\in C_c^\infty(\mathcal M)$, there exists
$C'''=C'''(f_1,f_2)>0$ such that whenever $t\geq1$ and $a_{tu+\sqrt t\,v}\in A_\te^+$,
$$
 \left|\Psi(t)e^{(2\rho_\te-\psi)(\log a_{tu+\sqrt t\,v})}\int_{\mathcal M}f_1(xa_{tu+\sqrt t\,v})f_2(x)\,dx\right|\leq C'''.
$$
\end{enumerate}
\end{theorem}
 
\begin{proof}
Set $a_t=a_{tu+\sqrt t\,v}$. 
By linearity, it suffices to treat nonnegative $f_1$ and $f_2$. 
Choose $\epsilon_0>0$ and $g_1,\ldots,g_\ell\in G$ such that each $\scrB_{\epsilon_0}(g_i)$ is contained in $\widetilde{\mathcal M}$, its projection to $\mathcal M$ is injective, and the sets $\Ga\scrB_{\epsilon_0/2}(g_i)$ cover $\supp f_1\cup\supp f_2$. Using a fixed smooth partition of unity, write
$$
 f_1=\sum_i f_{1,i},\qquad f_2=\sum_j f_{2,j},\qquad \supp f_{r,i}\subset\Ga\scrB_{\epsilon_0/2}(g_i).
$$
For every $k$, this gives
\begin{equation}\label{eq.fij}
 \sum_i\|f_{1,i}\|_{C^k}\ll\|f_1\|_{C^k},\qquad \sum_j\|f_{2,j}\|_{C^k}\ll\|f_2\|_{C^k}.
\end{equation}
The constants depend only on the fixed partition of unity.
 
Let $\tilde f_{2,j}\in C_c^\infty(\scrB_{\epsilon_0}(g_j))$ be the lift of $f_{2,j}$, extended by zero outside this set. Thus
$$
 f_{2,j}([g])=\sum_{\gamma\in\Ga}\tilde f_{2,j}(\gamma g.o).
$$
Set
$$
 D=(N_\te^-)_{\epsilon_0/2}(A_\te)_{\epsilon_0/2}.
$$
For $ha\in D$, define $\phi_{j,ha}(n)=\tilde f_{2,j}(g_jhan.o)$.
Note that $\phi_{j,ha}$ $(ha\in D)$ are supported on $\overline{(N_\te^+)_{\epsilon_0/2}}$ and 
$$
 \|\phi_{j,ha}\|_{C^k(N_\te^+)} \ll\|f_{2,j}\|_{C^k(\mathcal M)}.
$$
In these coordinates the Haar measure is $e^{-2\rho_\te(\log a)}\,dn\,da\,dh$. 
Hence
\begin{align}\label{eq.pu1}
 &\int_{\mathcal M}f_1(xa_t)f_2(x)\,dx=\sum_{i,j}\int_D\left(\int_{N_\te^+}f_{1,i}([g_jhan]a_t)\phi_{j,ha}(n)\,dn\right)e^{-2\rho_\te(\log a)}\,da\,dh.
\end{align}
 
Choose $\chi\in C_c^\infty((N_\te^+)_{\epsilon_0})$ with $0\leq\chi\leq1$ and $\chi=1$ on the common supports of $\phi_{j,ha}$. 
We first prove (2), assuming condition (2) of Definition \ref{lmp}. The uniformity in Proposition \ref{p2}(2), for $g_jha$ in the compact set
$\bigcup_jg_j\overline D$, gives constants $C_i$ such that, whenever $a_t\in A_\te^+$,
$$
 \Psi(t)e^{(2\rho_\te-\psi)(\log a_t)} \int_{N_\te^+}f_{1,i}([g_jhan]a_t)\chi(n)\,dn\leq C_i.
$$
Since $\phi_{j,ha}\leq\|f_{2,j}\|_\infty\chi$, it follows that
\begin{align*}
 &\Psi(t)e^{(2\rho_\te-\psi)(\log a_t)} \int_{\mathcal M}f_1(xa_t)f_2(x)\,dx \\ &\leq \left(\int_D e^{-2\rho_\te(\log a)}\,da\,dh\right) \left(\sum_iC_i\right)\left(\sum_j\|f_{2,j}\|_\infty\right).
\end{align*}
This proves (2).
 
For (1), assume only condition (1) of Definition \ref{lmp} and fix $v$. 
Use the upper bounds in the proofs of Propositions \ref{p1}(2) and \ref{p2}(2), with $\epsilon$ and $\delta_0$ fixed. 
As $g_jha$ ranges over $\bigcup_jg_j\overline D$, the functions used there can be chosen with common compact supports and uniformly bounded sup norms, since the denominators $\mu_{g,N_\te^\pm}^{\PS}(\rho)$ are continuous in $g$ and can be made uniformly positive by choosing $\rho$ on a finite cover. 
Hence there are nonnegative $B_i,\Phi\in C_c^\infty(\mathcal M)$ and $C>0$, independent of $ha$ and $t$, such that, for all sufficiently large $t$,
\begin{align*}
 &\Psi(t)e^{(2\rho_\te-\psi)(\log a_t)}\int_{N_\te^+}f_{1,i}([g_jhan]a_t)\chi(n)\,dn\\
 &\quad\leq C\Psi(t)\int_{\mathcal M}B_i(xa_t)\Phi(x)\,dm^{\BMS}(x).
\end{align*}
 Condition (1) of Definition \ref{lmp} bounds the right-hand side for this fixed $v$. Since $\phi_{j,ha}\leq\|f_{2,j}\|_\infty\chi$, the integrand is bounded by an integrable function on $D$.
 
Apply Proposition \ref{p2}(1) with $g$ replaced by $g_jha$ in \eqref{eq.pu1}. Dominated convergence and \eqref{BRPNdual} give
\begin{align}\label{eq.pou}
 &\lim_{t\to\infty}\Psi(t)e^{(2\rho_\te-\psi)(\log a_t)} \int_{\mathcal M}f_1(xa_t)f_2(x)\,dx\\
 &=J(v)\sum_i m^{\BR}(f_{1,i}) \sum_j\int_D \mu_{g_jha,N_\te^+}^{\PS}(\phi_{j,ha}) e^{-2\rho_\te(\log a)}\,da\,dh\notag\\
 &=J(v)\sum_i m^{\BR}(f_{1,i})\sum_jm^{\BR_*}(f_{2,j})=J(v)m^{\BR}(f_1)m^{\BR_*}(f_2).\notag
\end{align}
This proves (1).
\end{proof}
 By uniform approximation, (1) and (2) also hold for $f_1,f_2\in C_c(\mathcal M)$. The approximation errors are bounded by fixed nonnegative smooth functions, and the corresponding main terms are bounded by local finiteness of the BR measures. In particular, (1) proves Theorem \ref{thm.haar}.
\begin{theorem}
 Suppose that $\Gamma$ is linearly $\theta$-regular, and let $\mathcal K\subset\mathcal M$ be a compact set. Then (1) and (2) below hold under conditions (1) and (2) of Definition \ref{def.em}, respectively, for some $d'''>0$ depending only on $d$, $k$ and $G$.
\begin{enumerate}
\item For every $v\in\ker\psi$ and $f_1,f_2\in C_c^\infty(\cal K)$,
$$
 \err'''(a_{tu+\sqrt t\,v},f_1,f_2)=O_{\cal K,v}\left(t^{-d'''}\norm{f_1}_{C^k(\mathcal M)}\norm{f_2}_{C^k(\mathcal M)}\right)
 \qquad\text{as $t\to\infty$}.
$$
\item There exists $Q_{\cal K}'''\in(L^1\cap L^\infty)(\ker\psi)$ such that, for all $f_1,f_2\in C_c^\infty(\cal K)$, $t\geq1$ and $v\in\ker\psi$ with $a_{tu+\sqrt t\,v}\in A_\te^+$,
\begin{equation*}
 |\err'''(a_{tu+\sqrt t\,v},f_1,f_2)|\leq Q_{\cal K}'''(v)t^{-d'''}\norm{f_1}_{C^k(\mathcal M)}\norm{f_2}_{C^k(\mathcal M)}.
\end{equation*}
\end{enumerate}
\end{theorem}
\begin{proof}
Use the same partition of unity and notation as in the proof of Theorem \ref{prop.mixH0}, with $a_t=a_{tu+\sqrt t\,v}$.
We use the error term $\err''$ of Proposition \ref{p2}, with $g$ replaced by $g_jha$.
Multiplying \eqref{eq.pu1} by $\Psi(t)e^{(2\rho_\te-\psi)(\log a_t)}$ and subtracting the main term in \eqref{eq.pou} gives
\begin{equation*}
 \err'''(a_t,f_1,f_2) =\sum_{i,j}\int_D \err''(a_t,f_{1,i},\phi_{j,ha}) e^{-2\rho_\te(\log a)}\,da\,dh.
\end{equation*}
For (1), assume condition (1) of Definition \ref{def.em} and fix $v$. Proposition \ref{p2}(3) and \eqref{eq.fij} give
\begin{align*}
 |\err'''(a_t,f_1,f_2)| &\ll_{\mathcal K,v} t^{-d''}\sum_{i,j} \|f_{1,i}\|_{C^k}\|f_{2,j}\|_{C^k} \int_D e^{-2\rho_\te(\log a)}\,da\,dh \\ & \ll_{\mathcal K,v} t^{-d''}\|f_1\|_{C^k}\|f_2\|_{C^k}.
\end{align*}
Thus (1) holds with $d'''=d''$. For (2), assume condition (2) of Definition \ref{def.em}. The uniformity in Proposition \ref{p2}(4) gives one function
$Q''\in(L^1\cap L^\infty)(\ker\psi)$ for all $g_jha$ above. For $t\geq1$ and $a_t\in A_\te^+$, the
same calculation yields
$$
 |\err'''(a_t,f_1,f_2)| \leq C t^{-d''}Q''(v)\|f_1\|_{C^k}\|f_2\|_{C^k}.
$$
Taking $Q'''_{\mathcal K}=C Q''$ proves (2).  The constants in (1), and how large $t$ must be, do not depend on $f_1,f_2\in C_c^\infty(\mathcal K)$. Hence (1) also proves Theorem \ref{thm.eh}.
\end{proof}
The first part of the following corollary follows by approximating the characteristic functions of $E$ and $F$ by continuous functions. The second follows first for product functions by dividing by $\operatorname{vol}(R_t)$, and then for general $\Phi$ by uniform approximation on $\overline E\times\overline F$.
 
\begin{corollary}\label{cor.vol2}\rm Suppose that $m_\psi^{\BMS}$ satisfies local mixing with the pair $(\Psi,J)$.
Let $E,F\subset\mathcal M$ be relatively compact Borel sets such that $ m_\psi^{\BR_*}(\partial E)=m_\psi^{\BR}(\partial F)=0$ and $m_\psi^{\BR_*}(E)m_\psi^{\BR}(F)>0$.
For $t>0$, set
$
 R_t:=\{x\in E:xa_{tu}\in F\}=E\cap Fa_{-tu}.
$
Then as $t\to\infty$,
\begin{equation*}
 \operatorname{vol}(R_t)\sim \frac{J(0)m_\psi^{\BR_*}(E)m_\psi^{\BR}(F)}{\Psi(t)}e^{t(\psi_\Gamma^\theta(u)-2\rho_\theta(u))}.
\end{equation*}
Moreover, for every $\Phi\in C(\overline E\times\overline F)$,
\begin{equation*}
\lim_{t\to\infty}\frac{\int_{R_t}\Phi(x,xa_{tu})\,dx}{\operatorname{vol}(R_t)}
 =\frac{\int_E\int_F\Phi(x,y)\,dm_\psi^{\BR}(y)\,dm_\psi^{\BR_*}(x)}{m_\psi^{\BR_*}(E)m_\psi^{\BR}(F)}.
\end{equation*}
\end{corollary}

\begin{proof}[Proof of Corollary~\ref{cor.vol}]
Fix $p\in\mathcal M$ and $0\ne H\in\mathfrak a_\theta$, and put $X(x):=\left.\frac{d}{dt}\right|_{t=0}xa_{tH}$. For the squared distance $f(x):=d_{\mathcal M} (p,x)^2$, Riemannian normal coordinates at $p$ give $X^2f(p)=2\|X(p)\|^2>0$.
By continuity, choose a relatively compact product neighborhood of $p$ in local Hopf coordinates with $X=\partial_\tau$ and $\partial_\tau^2f>0$. For sufficiently small $r>0$, the closed ball $\overline{B(p,r)}$ lies in this neighborhood. With all other coordinates fixed, strict convexity implies that $f=r^2$ has at most two solutions in $\tau$. Since both BR measures are absolutely continuous in the $\tau$-coordinate, Fubini's theorem gives
$ m_\psi^{\BR_*}(\partial B(p,r))=m_\psi^{\BR}(\partial B(p,r))=0$. The support of $m_\psi^{\BMS}$ is contained in the supports of both BR measures. Thus sufficiently small balls centered at $p,q\in\operatorname{supp}m_\psi^{\BMS}$ satisfy the positivity and boundary conditions of  Corollary~\ref{cor.vol2}. Applying it with $E=B(p,r)$ and $F=B(q,s)$ proves the first claim.
 
For the second claim, put $C=J(0)m_\psi^{\BR_*}(E)m_\psi^{\BR}(F)>0$. Then
$
 \log\operatorname{vol}(R_t)=t\bigl(\psi_\Gamma^\theta(u)-2\rho_\theta(u)\bigr)-\log\Psi(t)+\log C+o(1).
$ Dividing by $t$ and using $\log\Psi(t)=o(t)$ (\Cref {lem.ratio}) proves \eqref{eq.growth-volume}. Finally, \eqref{avol} follows from \eqref{eq.rv} and Proposition \ref{prop.app}.
\end{proof}

\appendix
\section{Local mixing for general Anosov subgroups}\label{sec.app}
 
 In this appendix, we assume that $\Gamma$ is $\theta$-Anosov. We follow the proof of Chow-Sarkar \cite{CS}, using their notation, and indicate the changes needed for the partial Weyl chamber action. This gives the local mixing statement and the uniform bound used in this paper.
 
Fix $u\in\operatorname{int}\mathcal L_\theta$ and its tangent form $\psi$, and let $\nu=\nu_\psi$, $\nu_{\mathrm i}=\nu_{\psi\circ\i}$ and $m=m^{\BMS}$ be as in Section \ref{sec.set}. 
Let $ W=\ker\psi$, $q=\dim\mathfrak a_\theta$ and $ r=\frac{q-1}{2}$.
Following \cite[Section 3.2]{CS}, write
\begin{equation*} \check\Omega:=\Omega=\Gamma\backslash(\Lambda_\theta^{(2)}\times\mathfrak a_\theta),\qquad \mathcal X:=\Gamma\backslash(\Lambda_\theta^{(2)}\times\mathbb R),
\end{equation*}
where the action of $\Gamma$ on $\Lambda_\theta^{(2)}\times\mathbb R$ is given by $\gamma(\xi_1,\xi_2,s)=(\gamma\xi_1,\gamma\xi_2,s+\psi(\sigma_\theta(\gamma,\xi_1)))$.
The projection $\pi_\psi:\Omega\to \mathcal X$ and the translation flow on $\mathcal X$ are
\begin{equation*}
 \pi_\psi[\xi_1,\xi_2,b]=[\xi_1,\xi_2,\psi(b)],\qquad\phi_t[\xi_1,\xi_2,s]=[\xi_1,\xi_2,s+t].
\end{equation*}
In the arguments of \cite{CS}, replace $\mathfrak a,A,\sigma,\beta,\mu,\lambda$ by $\mathfrak a_\theta,A_\theta,\sigma_\theta,\beta^\theta,\mu_\theta,\lambda_\theta$, where $\lambda_\theta=p_\theta\circ\lambda$ and $\lambda(\gamma)=\lim_{n\to\infty}\mu(\gamma^n)/n$. 
 
\begin{proposition}\label{prop.app}\rm
There are $\kappa_u>0$ and a positive definite quadratic form $Q_u$ on $W$ such that, for every $v\in W$ and for $f_1,f_2\in C_c(\check\Omega)$,
\begin{equation}\label{eq.mix}
 t^r\int_{\check\Omega}f_1(xa_{tu+\sqrt t\,v})f_2(x)\,dm(x)\to \kappa_u e^{-Q_u(v)}m(f_1)m(f_2).
\end{equation}
The convergence is uniform for $v$ in compact subsets of $W$. Moreover, for each $f_1,f_2\in C_c(\check\Omega)$,
\begin{equation}\label{app.mixing-bound}
 \sup_{\substack{t>0\\v\in W}}t^r\left|\int_{\check\Omega}f_1(xa_{tu+\sqrt t\,v})f_2(x)\,dm(x)\right|<\infty.
\end{equation} Hence both conditions of Definition \ref{lmp} hold with $\Psi(t)=t^r$ and $J(v)=\kappa_u e^{-Q_u(v)}$. In particular, $J\in (L^1\cap L^\infty)(W)$.
\end{proposition}
 
\begin{proof}
We may assume that $\psi(u)=1$.
 
\noindent\textbf{The compact quotient and its Markov section.}
Since $\Gamma$ is $\theta$-Anosov, $\psi$ is positive on $\mathcal L_\theta\setminus\{0\}$; see \cite[Proposition 5.8.1 and Corollary 5.9.1]{Samb4}. By \cite[Corollary 5.3.3]{Samb4} (see also \cite[Theorem 9.2]{KOW}), $\mathcal X$ is compact and its translation flow is H\"older conjugate to a H\"older change of time of the Gromov geodesic flow.  Hence \cite[Theorems 3.5 and 4.15]{CS} hold in this setting.
The map $\pi_\psi$ is a bundle with fiber $W$, as in \cite[Section 3.2]{CS}. Let $B=\exp W$ and $\cal N:=\cal M/B$, which is a smooth manifold by Theorem \ref{theorem.dd}. Then $\cal X=\Omega/B$ is a closed subset of $\cal N$, invariant under the quotient flow.
 We use the restriction to $\check\Omega$ of the Riemannian distance fixed in Section \ref{sec.lcd}. On $\mathcal X$, we use the restriction of a Riemannian distance on the smooth manifold $\mathcal N$. On compact coordinate neighborhoods, these distances are bi-Lipschitz equivalent to the corresponding product distances in Hopf coordinates.
 
Writing $b=su+z$ with $z\in W$, the Hopf formula gives
\begin{equation*}
 dm=dm_0\,dz,
\end{equation*}
 where $m_0$ is the finite Hopf measure on $\mathcal X$, normalized by this equality. This is the product formula in \cite[Proposition 3.5]{Samb1} and \cite[Propositions 3.5.2-3.5.3]{Samb4}. The normalized measure $m_0/m_0(\mathcal X)$ is the measure of maximal entropy, and that entropy is $1$; see \cite[Proposition 3.3.2 and Remark 5.5.4]{Samb4}.
 
Define $W^{su}$ and $W^{ss}$ by \cite[(8)]{CS}, using $N_\theta^\pm$ and the partial Busemann maps. Their independence from the representative follows as in \cite[Lemma 4.16]{CS}, using \eqref{eq.basic1}. To obtain the Markov section of \cite[Corollary 4.17]{CS}, choose a Pl\"ucker representation associated to $P_\theta$. Its restriction to $\Gamma$ is projective Anosov \cite[Corollary 2.13]{BCLS}, 
so the associated flow is metric Anosov by \cite[Propositions 4.1 and 5.7]{BCLS};
see also \cite[Theorem 4.13]{CS}. Its weak stable and weak unstable leaves, like those defined from $W^{ss}$ and $W^{su}$, are determined by the first and second endpoints, respectively. The adjustment of rectangles along flow lines in the proof of Corollary 4.17 therefore gives a Markov section for these leaves.
 
Use $\mathcal R,\mathcal P,\Sigma,\Sigma^+,\sigma,\eta$ and $\tau$ as in \cite[Sections 4.1 and 4.3]{CS}. In particular, $\eta:\Sigma\to\mathcal R$ is the coding map, and $\tau$ is the H\"older extension of the first-return time from sequences for which $\eta(\sigma^nx)$ lies in the rectangle interiors for every $n\in\mathbb Z$. Let $p$ be the equilibrium probability measure for $-\tau$, and set $\bar\tau=\int_\Sigma\tau\,dp$.
 
\medskip
\noindent\textbf{The return vector and its periods.}
Choose lifts of the small Markov rectangles to $\Lambda_\theta^{(2)}\times\mathbb R$, and lift their points to $\check\Omega$ by $(\xi_1,\xi_2,s)\mapsto[\xi_1,\xi_2,su]$. Denote this section by $F:\mathcal R\to\check\Omega$. The return vector of \cite[Definition 5.4]{CS} is determined by
\begin{equation*}
 F(\mathcal Pz)=F(z)\exp\mathsf K(z),\qquad\psi(\mathsf K)=\tau,\qquad\widehat{\mathsf K}:=\mathsf K-\tau u\in W.
\end{equation*}
 
We check that $\mathsf K$ is constant on local strong stable leaves, as needed in \cite[Lemma 5.6]{CS}. Write the chosen lifts as $\tilde z=(\xi_1,\xi_2,s)$ and $\widetilde{\mathcal Pz}=(\xi'_1,\xi'_2,s')$, with $\widetilde{\mathcal Pz}=\gamma\phi_{\tau(z)}\tilde z$. 
Then $\mathsf K(z)=(s'-s)u-\sigma_\theta(\gamma,\xi_1).$
Along such a leaf, $\xi_1,s,\tau,\xi'_1,s'$ are fixed. They also determine $\gamma$: a nontrivial element fixing a limit flag has nonzero $\psi$-period, so the action on $\Lambda_\theta\times\mathbb R$ is free. Thus $\mathsf K$ and $\widehat{\mathsf K}$ are constant on these leaves and induce H\"older functions on $\Sigma^+$, depending only on the nonnegative coordinates.
 
The period identity in \cite[Lemma 6.6]{CS} becomes
\begin{equation}\label{eq.per}
 \mathsf K_n(x)=\lambda_\theta(\gamma),\qquad\tau_n(x)=\psi(\lambda_\theta(\gamma))\qquad(\sigma^nx=x),
\end{equation}
where the closed orbit and its multiplicity determine $\gamma$ as in that lemma. This is the partial Iwasawa period identity \cite[Lemma 5.3.2]{Samb4}.  We need the following density statement from \cite[Proposition 6.7]{CS}:
\begin{equation}\label{eq.per2}
 \overline{\left\langle\mathsf K_n(x):\sigma^nx=x,\ n\geq1\right\rangle}=\mathfrak a_\theta.
\end{equation}
 The proof of \cite[Proposition 6.7]{CS} applies after centering the first rectangle on a closed orbit corresponding to an element whose Jordan projection lies in $\operatorname{int}\mathfrak a^+$.  As in that proof, choose a Zariski dense Schottky subsemigroup whose full fixed flags are close to those of this element. Its partial fixed flags lie in the endpoint neighborhoods of the first rectangle, so its closed orbits are represented by periodic sequences through that rectangle. By \cite[Proposition 6.5]{CS}, its full Jordan projections generate a dense subgroup of $\mathfrak a$. Applying $p_\theta$ and using \eqref{eq.per} proves \eqref{eq.per2}.
 
\medskip
\noindent\textbf{The transfer operators and the mean return vector.}
Use the inner product $\langle\cdot,\cdot\rangle_\psi$ of \cite[Section 5]{CS}, with associated norm $\|\cdot\|_\psi$, so that $\psi(H)=\langle H,u\rangle_\psi$. The transfer operators are
\begin{equation*}
 (\mathcal L_\zeta H)(x)=\sum_{\sigma y=x}e^{-\tau(y)+i\langle\zeta,\mathsf K(y)\rangle_\psi}H(y),\qquad\zeta\in\mathfrak a_\theta.
\end{equation*}
After conjugating by the positive eigenfunction at $\zeta=0$, we have $\mathcal L_0 1=1$. The proofs of \cite[Proposition 6.9 and Theorems 6.10 and 6.13(2)]{CS} apply to these operators, using \eqref{eq.per2}. In particular, $1$ is not an eigenvalue for $\zeta\ne0$, and $(I-\mathcal L_\zeta)^{-1}$ is analytic away from $0$. Near $0$ it has the leading-eigenvalue decomposition of Theorem 6.13(2).
 We next show that the mean return vector is proportional to $u$. Set
$$
 M=\int_\Sigma\mathsf K\,dp,\qquad P(\varphi)=P_\sigma(-\varphi(\mathsf K)).
$$
Near $\psi$, the zero set $\mathcal C=\{P=0\}$ consists of the tangent forms; see \cite[Corollaries 2.6.4 and 5.9.1]{Samb4}. Differentiation of pressure gives
$
 DP_\psi(\omega)=-\omega(M)$ and $ DP_\psi(\psi)=-\bar\tau\ne0.
$
Thus $T_\psi\mathcal C=\{\omega:\omega(M)=0\}$. Every $\varphi\in\mathcal C$ satisfies $\varphi(u)\geq\psi_\Gamma^\theta(u)=\psi(u)=1$. Hence $\varphi\mapsto\varphi(u)$ has a minimum at $\psi$, and every $\omega\in T_\psi\mathcal C$ also satisfies $\omega(u)=0$. Both hyperplanes have codimension one, so $M$ is proportional to $u$. Applying $\psi$ gives
\begin{equation*}
 \int_\Sigma\mathsf K\,dp=\bar\tau u,\qquad\int_\Sigma\widehat{\mathsf K}\,dp=0.
\end{equation*}
The proof of \cite[Proposition 7.7]{CS} now applies. In particular, if $\kappa(\zeta)$ is the leading eigenvalue near $0$, there is a positive definite symmetric bilinear form $B$ on $\mathfrak a_\theta$ such that
$$ \kappa(\zeta)=1+i\bar\tau\langle\zeta,u\rangle_\psi-B(\zeta,\zeta)+o(\|\zeta\|_\psi^2).
$$
 
\medskip
\noindent\textbf{Local mixing.}
In the suspension of \cite[Section 7.3]{CS}, omit the $M_\Gamma$ coordinate. The measure on $0\leq s<\tau(x)$ is
$ \frac{m_0(\mathcal X)}{\bar\tau}\,dp(x)\,dz\,ds$.
The proof of \cite[Theorem 7.1]{CS} now applies with $\mathfrak a_\theta$ in place of $\mathfrak a$. 
The Fourier estimates and passage to continuous functions are given in \cite[Theorem 7.13(1)-(2), Lemma 7.14 and Proposition 7.15]{CS}. 
Taking the transverse scale to be $\sqrt t$ gives \eqref{eq.mix}, with a positive definite quadratic form $Q_u$ on $W$.  
The same Gaussian computation as in \cite[Remark B.4]{Samb4} gives the dependence on $v$, and the convergence is uniform on compact subsets of $W$ by \cite[Theorem 7.1(3)]{CS}.
 
The uniform estimates in the proof of \cite[Theorem 7.1(1)]{CS} give \eqref{app.mixing-bound} for large $t$, uniformly in $v\in W$. For bounded $t>0$, invariance of $m$ gives, for every $a\in\mathfrak a_\theta$,
$
 \left|\int_{\check\Omega}f_1(x\exp a)f_2(x)\,dm(x)\right|\leq\|f_1\|_{L^2(m)}\|f_2\|_{L^2(m)}.
$
This proves \eqref{app.mixing-bound}. Finally, restricting $C_c(\mathcal M)$ functions to the closed support $\check\Omega=\Omega$ gives both conditions in Definition \ref{lmp}, and positive definiteness of $Q_u$ gives $J\in (L^1\cap L^\infty)(W)$.
\end{proof}

\end{document}